\documentclass[12pt]{amsart}

\usepackage{amssymb}
\input amssym.def
\usepackage{pstricks}
\newtheorem{theorem}{Theorem}[section]

\newtheorem{prop}[theorem]{Proposition}

\theoremstyle{definition}
\newtheorem{rem}[theorem]{Remark}
\newcommand{\foorp}{\hfill$\Box$}
\newcommand{\M}{\mathcal{M}}
\newcommand{\Z}{\mathbb{Z}}
\newcommand{\R}{\mathbb{R}}
\newcommand{\T}{\mathcal{T}}
\newcommand{\Dif}{\mathcal{H}}
\newcommand{\bdr}[1]{\partial\! #1}
\newcommand{\C}{\mathcal{C}}
\newcommand{\Co}{\mathcal{C}}

\newcommand{\Stab}{\mathrm{Stab}}
\numberwithin{equation}{section}
\author{B{\l}a\.zej Szepietowski}
\title[Presentation for the mapping class group.]
{A presentation for the mapping class group of the closed
non-orientable surface of genus 4}
\address[]{Institute of Mathematics, Gda\'nsk University, Wita Stwosza 57,
80-952 Gda\'nsk, Poland} \email{blaszep@math.univ.gda.pl}
\begin{document}
%
%
\begin{abstract}
In \cite{Sz1} we proposed a method of finding a finite presentation
for the mapping class group of a non-orientable surface by using its
action on the so called ordered complex of curves. In this paper we
use this method to obtain an explicit finite presentation for the
mapping class group of the closed non-orientable surface of genus
$4$. The set of generators in this presentation consists of $5$ Dehn
twists, $3$ crosscap transpositions and one involution, and it can
be immediately reduced to the generating set found by Chillingworth
\cite{C}.
\end{abstract}

\maketitle
\section{Introduction}
Presentations for the mapping class group $\M(F)$ of an orientable
surface $F$ have been found by various authors. McCool \cite{McC}
was the first who showed that $\M(F)$ is finally presented. His
proof is purely algebraic and no concrete presentation was derived
from it. Hatcher and Thurston \cite{HT} showed how to obtain a
finite presentation for $\M(F)$ from its action on a simply
connected $2$-dimensional complex. Using their result, Wajnryb
\cite{W} obtained a simple presentation for $\M(F)$, for $F$
having at most one boundary component. Starting from Wajnryb's
result, Gervais \cite{G} found a finite presentation for $\M(F)$,
for $F$ having genus at least one and arbitrary many boundary
components. Benvenuti \cite{B} showed how the Gervais presentation
may be recovered by using the so called ordered complex of curves,
which is a modification of the classical complex of curves defined
by Harvey \cite{H}, instead of the complex of Hatcher and
Thurston. In \cite{Sz1} we used Benvenuti's approach to obtain a
presentation for the mapping class group of an arbitrary compact
non-orientable surface, defined in terms of mapping class group of
complementary surfaces of collections of simple closed curves.  In
this paper we find an explicit finite presentation for the mapping
class group of a closed non-orientable surface of genus $4$, by
using results of \cite{Sz1}. It is very difficult to derive an
explicit presentation for $\M(F)$ for general $F$ from the
presentation in \cite{Sz1} because of its recursive form. The
number of subsurfaces involved in the presentation increases with
the genus and the number of boundary components of $F$.
Furthermore, even if one is only interested in the case when $F$
is closed, one still has to consider surfaces with boundary
obtained by cutting, which appear to be more difficult to handle.

In contrast to the case of orientable surfaces, little is known
about the mapping class group $\M(F)$ of a non-orientable surface
$F$. In particular, no explicit finite presentation for $\M(F)$ is
known if $F$ has genus at least $4$. If $F$ is closed and has
genus $g$, then $M(F)$ is trivial if $g=1$ and isomorphic to
$\mathbb{Z}_2\times\mathbb{Z}_2$ if $g=2$ (see \cite{L}). For
$g=3$ a simple presentation for $\M(F)$ was found by Birman and
Chillingworth \cite{BC}. Lickorish \cite{L,L1} proved that $\M(F)$
is generated by Dehn twists and one crosscap slide (or
Y-homeomorphism) if $g\ge 2$ and Chillingworth \cite{C} found a
finite generating set for $\M(F)$. If $F$ is not closed, then a
finite set of generators for $\M(F)$ was found by Korkmaz \cite{K}
if $F$ has punctures, and by Stukow \cite{S0} if $F$ has punctures
and boundary and $g\ge 3$.

This paper is organized as follows. In the next section we present
basic definitions and facts and state our main result, Theorem
\ref{main}, which is a presentation for the mapping class group
$\M(F)$ of a closed non-orientable surface $F$ of genus $4$. We
also show that the proposed relations hold in $\M(F)$. In Section
\ref{s_pres} we determine orbits of the action of $\M(F)$ on the
ordered complex of curves $\C$ and describe a presentation for
$\M(F)$ arising from this action. In Section \ref{s_stab} we
determine stabilizers of vertices and edges of $\C$. Finally, in
Section \ref{inj} we show that relations in Theorem \ref{main} are
indeed defining relations for $\M(F)$.

\section{\label{preli}Preliminaries}
\subsection{Basic definitions.}
Let $F$ denote a connected surface, orientable or not, possibly with
boundary. Define $\Dif(F)$ to be the group of all (orientation
preserving if $F$ is orientable) homeomorphisms $h\colon F\to F$
equal to the identity on the boundary of $F$. The {\it mapping class
group} $\M(F)$ is the group of isotopy classes in $\Dif(F)$. By
abuse of notation we will use the same symbol to denote a
homeomorphism and its isotopy class. If $g$ and $h$ are two
homeomorphisms, then the composition $gh$ means that $h$ is applied
first. In this paper all surfaces and curves are assumed to have
PL-structure, and all homeomorphisms, embeddings and isotopies are
piecewise linear.

\medskip

By a {\it simple closed curve} in $F$ we mean an embedding
$\gamma\colon S^1\to F$. Note that $\gamma$ has an orientation; the
curve with opposite orientation but same image will be denoted by
$\gamma^{-1}$. By abuse of notation, we also use $\gamma$ for the
image of $\gamma$. If $\gamma_1$ and $\gamma_2$ are isotopic, we
write $\gamma_1\simeq\gamma_2$.

We say that $\gamma$ is {\it non-separating} if $F\backslash\gamma$
is connected and {\it separating} otherwise. According to whether a
regular neighborhood of $\gamma$ is an annulus or a M\"obius strip,
we call $\gamma$ respectively {\it two-} or {\it one-sided}. If
$\gamma$ is one-sided, then we denote by $\gamma^2$ its double, i.e.
the curve $\gamma^2(z)=\gamma(z^2)$ for $z\in S^1\subset\mathbb{C}$.
Note that although $\gamma^2$ is not simple, it is freely homotopic
to a two-sided simple closed curve.

We say that $\gamma$ is {\it generic} if it neither bounds a disk
nor a M\"obius strip. 

Define a {\it generic n-family of disjoint curves} to be an ordered
$n$-tuple $(\gamma_1,\dots,\gamma_n)$ of generic simple closed
curves satisfying:
\begin{itemize}
\item $\gamma_i\cap\gamma_j=\emptyset$, for $i\ne j$;
\item $\gamma_i$ is neither isotopic to $\gamma_j$ nor to $\gamma_j^{-1}$, for $i\ne j$.
\end{itemize}

We say that two generic $n$-families of disjoint curves
$(\gamma_1,\dots,\gamma_n)$ and $(\gamma'_1,\dots,\gamma'_n)$ are
{\it equivalent} if $\gamma_i\simeq(\gamma'_i)^{\pm 1}$ for each
$1\le i\le n$. We write $[\gamma_1,\dots,\gamma_n]$ for the
equivalence class of a generic $n$-family of disjoint curves.

The {\it ordered complex of curves} of $F$ is the $\Delta$-complex
(in the sens of \cite{Hat}, Chapter 2) whose $n$-simplices are the
equivalence classes of generic $(n+1)$-families of disjoint curves
in $F$. If $[\gamma_1,\dots,\gamma_{n+1}]$ is $n$-simplex then its
faces are the $(n-1)$-simplices
$[\gamma_1,\dots,\widehat{\gamma_i},\dots,\gamma_{n+1}]$ for
$i=1,\dots,n+1$, where $\widehat{\gamma_i}$ means that $\gamma_i$ is
deleted. We denote this complex by $\C$. Simplices of dimension $0$,
$1$ and $2$ are called vertices, edges and triangles respectively.
Vertices of $\C$ are the isotopy classes of unoriented generic
curves. The mapping class group $\M(F)$ acts on $\C$ by
$h[\gamma_1,\dots,\gamma_r]=[h(\gamma_1),\dots,h(\gamma_n)]$.

The idea of the ordered complex of curves comes from \cite{B}. It is
a variation of the classical complex of curves introduced by Harvey
\cite{H}. 

\medskip

Given a two-sided simple closed curve $\gamma$ we can define a Dehn
twist $c$ about $\gamma$. On a  non-orientable surface it is
impossible to distinguish between right and left twists, thus the
direction of a twist $c$ has to be specified for each curve
$\gamma$. Equivalently we may choose an orientation of a tubular
neighborhood of $\gamma$. Then $c$ denotes the right Dehn twist with
respect to the chosen orientation. Unless we specify which of the
two twists we mean, $c$ denotes (the isotopy class of) any of the
two possible twists.

\medskip

\begin{figure}
\input{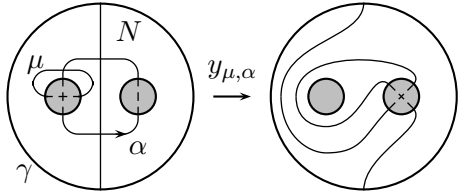}
\caption{\label{Y}Crosscap slide.}
\end{figure}

Suppose that $\mu$ and $\alpha$ are two simple closed curves in $F$,
such that $\mu$ is one-sided, $\alpha$ is two-sided and they
intersect at one point. Let $N$ be a regular neighborhood of
$\mu\cup\alpha$, which is homeomorphic to the Klein bottle with a
hole, and let $M\subset N$ be a regular neighborhood of $\mu$, which
is a M\"obius strip. We denote by $y_{\mu,\alpha}$ the {\it
Y-homeomorphism}, or {\it crosscap slide} of $N$ which may be
described as a result of sliding $M$ once along $\alpha$ keeping the
boundary of $N$ fixed. Figure \ref{Y} illustrates the effect of
$y_{\mu,\alpha}$ on an arc connecting two points in the boundary of
$N$. Here, and also in other figures of this paper, the shaded discs
represent crosscaps; this means that their interiors should be
removed, and then antipodal points in each resulting boundary
component should be identified. The homeomorphism $y_{\mu,\alpha}$
pushes the left crosscap through the right one, along $\alpha$.
Y-homeomorphism was first introduced by Lickorish; see \cite{L} for
a formal definition. Observe that $y_{\mu,\alpha}$ reverses the
orientation of $\mu$.  We extend $y_{\mu,\alpha}$ by the identity
outside $N$ to a homoeomorphism of $F$, which we denote by the same
symbol. Up to isotopy, $y_{\mu,\alpha}$ does not depend on the
choice of $N$. It also does not depend on the orientation of $\mu$
but does depend on the orientation of $\alpha$. The following
properties of Y-homeomorphisms are easy to verify:
\begin{equation}\label{y^-1}
y_{\mu,\alpha^{-1}}=y^{-1}_{\mu,\alpha};
\end{equation}
\begin{equation}\label{y2}
y_{\mu,\alpha}^2=c,
\end{equation}
where $c$ is Dehn twist about $\gamma=\bdr{N}$, right with respect
to the standard orientation of the plane of Figure \ref{Y};
\begin{equation}\label{hyh}
hy_{\mu,\alpha}h^{-1}=y_{h(\mu),h(\alpha)},
\end{equation}
for all $h\in\Dif(F)$.

\begin{figure}
\input{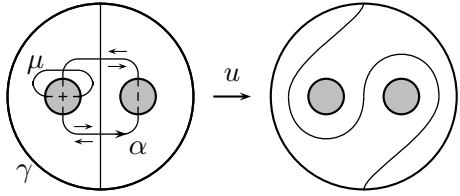}
\caption{\label{U}Crosscap transposition.}
\end{figure}

Let $a$ denote Dehn twist about $\alpha$ in direction indicated by
arrows in Figure \ref{U}. Then $u=ay_{\mu,\alpha}$ interchanges two
crosscaps (Figure \ref{U}). We call this homeomorphism {\it crosscap
transposition}. Since $u$ reverses orientation of a neighborhood of
$\alpha$, thus
\begin{equation}\label{uau}
uau^{-1}=a^{-1},
\end{equation}
\begin{equation}\label{u2}
u^2=y_{\mu,\alpha}^2=c.
\end{equation}

\subsection{Relations in $\M(F)$}
\begin{figure}
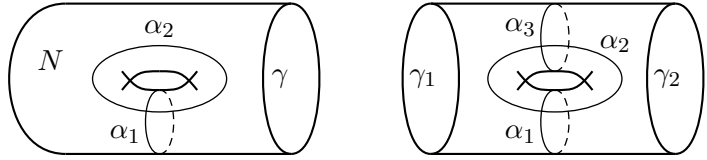

\begin{tabular}{cc}
\input{torus1}&\input{torus2}
\end{tabular}
\caption{\label{torus}Torus with one and two holes.}
\end{figure}
Suppose that $\alpha_1$ and $\alpha_2$ are two-sided simple closed
curves in $F$, intersecting at one point. Let $N$ be oriented
regular neighborhood of $\alpha_1\cup\alpha_2$, which is torus with
a hole, and let $\gamma$ denote its boundary (Figure \ref{torus}).
If $a_1$, $a_2$ and $c$ are Dehn twist about $\alpha_1$, $\alpha_2$
and $\gamma$ respectively, right with respect to the orientation of
$N$, then the following relations hold in $\M(F)$:
\begin{equation}\label{braid}
a_1a_2a_1=a_2a_1a_2,
\end{equation}
\begin{equation}\label{star1}
(a_1^2a_2)^4=c.
\end{equation}
First is the well known ``braid'' relation, second is a special case
of the ``star'' relation (see \cite{G}).

Consider the torus with two holes in the right hand side of Figure
\ref{torus} as embedded in $F$. If $a_1$, $a_2$, $a_3$, $c_1$ and
$c_2$ are Dehn twists about $\alpha_1$, $\alpha_2$, $\alpha_3$,
$\gamma_1$ and $\gamma_2$ respectively, right with respect to some
orientation of the torus, then the following relation holds in
$\M(F)$:
\begin{equation}\label{star2}
(a_1a_2a_3)^4=c_1c_2.
\end{equation}
This is also a special case of the ``star'' relation.

\medskip

\begin{figure}
\begin{tabular}{cc}
\pspicture*(3.75,3.2)
\pscircle(2,1.6){1.5}\rput(.65,.45){\small$\gamma_1$}
\pscircle(2,2.4){.3}\rput(2,2.4){\small$\gamma_2$}
\pscircle*[linecolor=lightgray](1.4,1.2){.3} \pscircle(1.4,1.2){.3}
\pscircle*[linecolor=lightgray](2.6,1.2){.3} \pscircle(2.6,1.2){.3}
\psline[linewidth=.5pt, arrowsize=2pt 2.5]{->}(2.15,.65)(2.3,.65)
\psline[linewidth=.5pt,linearc=.2](1.4,1.5)(1.4,1.8)(2.6,1.8)(2.6,1.5)
\psline[linewidth=.5pt,linearc=.2](1.4,.9)(1.4,.65)(2.6,.65)(2.6,.9)
\psline[linewidth=.5pt,linestyle=dashed,dash=3pt
2pt](1.4,.9)(1.4,1.5)
\psline[linewidth=.5pt,linestyle=dashed,dash=3pt
2pt](2.6,.9)(2.6,1.5) \rput[t](1.95,.55){\small$\alpha_1$}
\psline[linewidth=.4pt]{<-}(1.9,1.9)(2.2,1.9)
\psline[linewidth=.4pt]{->}(1.9,1.7)(2.2,1.7)
\psline[linewidth=.5pt,linestyle=dashed,dash=3pt
2pt](1.1,1.2)(1.7,1.2)
\psline[linewidth=.5pt,linearc=.2](1.1,1.2)(.95,1.2)(.95,1.6)(1.85,1.6)(1.85,1.2)(1.7,1.2)
\rput[br](1,1.55){\small$\mu$}
%
%
%
%
\endpspicture&\input{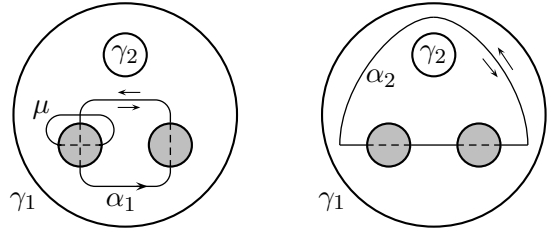}
\end{tabular}
\caption{\label{K}Klein bottle with two holes.}
\end{figure}

Consider the Klein bottle with two holes in Figure \ref{K} as
embedded in $F$. Let $a_1$ and $a_2$ denote Dehn twists about
$\alpha_1$ and $\alpha_2$ respectively, in the indicated directions.
Let $c_1$, $c_2$ denote Dehn twists about $\gamma_1$, $\gamma_2$,
right with respect to the standard orientation of the plane of the
figure and let $u$ denote the crosscap transposition
$u=a_1y_{\mu,\alpha_1}$. Then, by Lemma 7.8 in \cite{Sz1}, the
following relation holds in $\M(F)$:
\begin{equation}\label{Krel}
(ua_2)^2=c_1c_2.
\end{equation}

\subsection{Statement of the main result.}

\begin{figure}
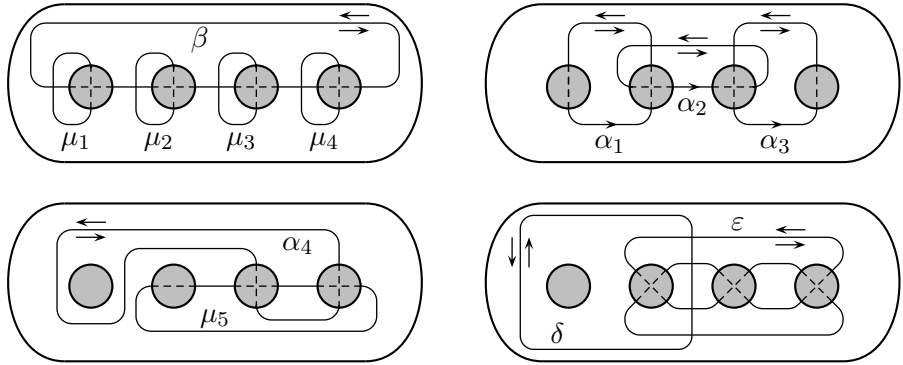

\begin{tabular}{cc}
\pspicture*(6,2.5)
\psline(1,.25)(5,.25) \psline(1,2.35)(5,2.35)
\psbezier(1,.25)(0,.25)(0,2.35)(1,2.35)
\psbezier(5,.25)(6,.25)(6,2.35)(5,2.35)
\pscircle*[linecolor=lightgray](1.35,1.25){.3}
\pscircle(1.35,1.25){.3}
\pscircle*[linecolor=lightgray](2.45,1.25){.3}
\pscircle(2.45,1.25){.3}
\pscircle*[linecolor=lightgray](3.55,1.25){.3}
\pscircle(3.55,1.25){.3}
\pscircle*[linecolor=lightgray](4.65,1.25){.3}
\pscircle(4.65,1.25){.3}
\rput[t](2.8,2.05){\small$\beta$}
\psline[linewidth=.5pt,linearc=.2](1.05,1.25)(.55,1.25)(.55,2.1)(5.45,2.1)(5.45,1.25)(4.95,1.25)
\psline[linewidth=.5pt,linestyle=dashed,dash=3pt
2pt](1.05,1.25)(1.65,1.25)
\psline[linewidth=.5pt](1.65,1.25)(2.15,1.25)
\psline[linewidth=.5pt,linestyle=dashed,dash=3pt
2pt](2.15,1.25)(2.75,1.25)
\psline[linewidth=.5pt](2.75,1.25)(3.25,1.25)
\psline[linewidth=.5pt,linestyle=dashed,dash=3pt
2pt](3.25,1.25)(3.85,1.25)
\psline[linewidth=.5pt](3.85,1.25)(4.35,1.25)
\psline[linewidth=.5pt,linestyle=dashed,dash=3pt
2pt](4.35,1.25)(4.95,1.25)
\psline[linewidth=.5pt, arrowsize=2pt 2.5]{<-}(4.65,2.2)(5.05,2.2)
\psline[linewidth=.5pt, arrowsize=2pt 2.5]{->}(4.65,2)(5.05,2)
\rput[t](1.15,.65){\small$\mu_1$}
\psline[linewidth=.5pt,linearc=.2](1.35,1.55)(1.35,1.7)(.85,1.7)(.85,.75)(1.35,.75)(1.35,.95)
\psline[linewidth=.5pt,linestyle=dashed,dash=3pt
2pt](1.35,.95)(1.35,1.55)
\rput[t](2.25,.65){\small$\mu_2$}
\psline[linewidth=.5pt,linearc=.2](2.45,1.55)(2.45,1.7)(1.95,1.7)(1.95,.75)(2.45,.75)(2.45,.95)
\psline[linewidth=.5pt,linestyle=dashed,dash=3pt
2pt](2.45,.95)(2.45,1.55)
\rput[t](3.35,.65){\small$\mu_3$}
\psline[linewidth=.5pt,linearc=.2](3.55,1.55)(3.55,1.7)(3.05,1.7)(3.05,.75)(3.55,.75)(3.55,.95)
\psline[linewidth=.5pt,linestyle=dashed,dash=3pt
2pt](3.55,.95)(3.55,1.55)
\rput[t](4.45,.65){\small$\mu_4$}
\psline[linewidth=.5pt,linearc=.2](4.65,1.55)(4.65,1.7)(4.15,1.7)(4.15,.75)(4.65,.75)(4.65,.95)
\psline[linewidth=.5pt,linestyle=dashed,dash=3pt
2pt](4.65,.95)(4.65,1.55)
%

%
\endpspicture&\pspicture*(6,2.5)
\psline(1,.25)(5,.25) \psline(1,2.35)(5,2.35)
\psbezier(1,.25)(0,.25)(0,2.35)(1,2.35)
\psbezier(5,.25)(6,.25)(6,2.35)(5,2.35)
\pscircle*[linecolor=lightgray](1.35,1.25){.3}
\pscircle(1.35,1.25){.3}
\pscircle*[linecolor=lightgray](2.45,1.25){.3}
\pscircle(2.45,1.25){.3}
\pscircle*[linecolor=lightgray](3.55,1.25){.3}
\pscircle(3.55,1.25){.3}
\pscircle*[linecolor=lightgray](4.65,1.25){.3}
\pscircle(4.65,1.25){.3}
\psline[linewidth=.5pt, arrowsize=2pt 2.5]{->}(1.9,.75)(2,.75)
\rput[t](1.9,.6){\small$\alpha_1$}
\psline[linewidth=.5pt,linearc=.2](1.35,1.55)(1.35,2.1)(2.45,2.1)(2.45,1.55)
\psline[linewidth=.5pt,linearc=.2](1.35,.95)(1.35,.75)(2.45,.75)(2.45,.95)
\psline[linewidth=.5pt,linestyle=dashed,dash=3pt
2pt](1.35,.95)(1.35,1.55)
\psline[linewidth=.5pt, arrowsize=2pt 2.5]{<-}(1.7,2.2)(2.1,2.2)
\psline[linewidth=.5pt, arrowsize=2pt 2.5]{->}(1.7,2)(2.1,2)
\psline[linewidth=.5pt, arrowsize=2pt 2.5]{->}(4.1,.75)(4.2,.75)
\rput[t](4.1,.6){\small$\alpha_3$}
\psline[linewidth=.5pt,linearc=.2](3.55,1.55)(3.55,2.1)(4.65,2.1)(4.65,1.55)
\psline[linewidth=.5pt,linearc=.2](3.55,.95)(3.55,.75)(4.65,.75)(4.65,.95)
\psline[linewidth=.5pt,linestyle=dashed,dash=3pt
2pt](1.35,.95)(1.35,1.55)
\psline[linewidth=.5pt,linestyle=dashed,dash=3pt
2pt](2.45,.95)(2.45,1.55)
\psline[linewidth=.5pt,linestyle=dashed,dash=3pt
2pt](3.55,.95)(3.55,1.55)
\psline[linewidth=.5pt, arrowsize=2pt 2.5]{<-}(3.9,2.2)(4.3,2.2)
\psline[linewidth=.5pt, arrowsize=2pt 2.5]{->}(3.9,2)(4.3,2)
\psline[linewidth=.5pt, arrowsize=2pt 2.5]{->}(3,1.25)(3.1,1.25)
\rput[t](3,1.1){\small$\alpha_2$}
\psline[linewidth=.5pt,linearc=.2](2.15,1.25)(2,1.25)(2,1.8)(4,1.8)(4,1.25)(3.85,1.25)
\psline[linewidth=.5pt,linestyle=dashed,dash=3pt
2pt](2.15,1.25)(2.75,1.25)
\psline[linewidth=.5pt](2.75,1.25)(3.25,1.25)
\psline[linewidth=.5pt,linestyle=dashed,dash=3pt
2pt](3.25,1.25)(3.85,1.25)
\psline[linewidth=.5pt, arrowsize=2pt 2.5]{<-}(2.8,1.9)(3.2,1.9)
\psline[linewidth=.5pt, arrowsize=2pt 2.5]{->}(2.8,1.7)(3.2,1.7)
\psline[linewidth=.5pt,linestyle=dashed,dash=3pt
2pt](4.65,.95)(4.65,1.55)
%

%
\endpspicture\\
\pspicture*(6,2.5)
\psline(1,.25)(5,.25) \psline(1,2.35)(5,2.35)
\psbezier(1,.25)(0,.25)(0,2.35)(1,2.35)
\psbezier(5,.25)(6,.25)(6,2.35)(5,2.35)
\pscircle*[linecolor=lightgray](1.35,1.25){.3}
\pscircle(1.35,1.25){.3}
\pscircle*[linecolor=lightgray](2.45,1.25){.3}
\pscircle(2.45,1.25){.3}
\pscircle*[linecolor=lightgray](3.55,1.25){.3}
\pscircle(3.55,1.25){.3}
\pscircle*[linecolor=lightgray](4.65,1.25){.3}
\pscircle(4.65,1.25){.3}
\rput[b](3,.7){\small$\mu_5$}
\psline[linewidth=.5pt,linearc=.2](2.15,1.25)(1.95,1.25)(1.95,.65)(5.15,.65)(5.15,1.25)(4.95,1.25)
\psline[linewidth=.5pt,linestyle=dashed,dash=3pt
2pt](2.15,1.25)(2.75,1.25)
\psline[linewidth=.5pt](2.75,1.25)(3.25,1.25)
\psline[linewidth=.5pt,linestyle=dashed,dash=3pt
2pt](3.25,1.25)(3.85,1.25)
\psline[linewidth=.5pt](3.85,1.25)(4.35,1.25)
\psline[linewidth=.5pt,linestyle=dashed,dash=3pt
2pt](4.35,1.25)(4.95,1.25)
\psline[linewidth=.5pt,linestyle=dashed,dash=3pt
2pt](3.55,.95)(3.55,1.55)
\psline[linewidth=.5pt,linestyle=dashed,dash=3pt
2pt](4.65,.95)(4.65,1.55)
\psline[linewidth=.5pt,linearc=.2](3.55,.95)(3.55,.8)(4.65,.8)(4.65,.95)
\psline[linewidth=.5pt,linearc=.2](3.55,1.55)(3.55,1.75)(1.8,1.75)(1.8,.75)(.9,.75)(.9,2)(4.65,2)(4.65,1.55)
\rput[t](4.1,1.9){\small$\alpha_4$}
\psline[linewidth=.5pt, arrowsize=2pt 2.5]{<-}(1.15,2.1)(1.55,2.1)
\psline[linewidth=.5pt, arrowsize=2pt 2.5]{->}(1.15,1.9)(1.55,1.9)
%

%
\endpspicture&\pspicture*(6,2.5)
\psline(1,.25)(5,.25) \psline(1,2.35)(5,2.35)
\psbezier(1,.25)(0,.25)(0,2.35)(1,2.35)
\psbezier(5,.25)(6,.25)(6,2.35)(5,2.35) %
\pscircle*[linecolor=lightgray](1.35,1.25){.3}
\pscircle(1.35,1.25){.3}
\pscircle*[linecolor=lightgray](2.45,1.25){.3}
\pscircle(2.45,1.25){.3}
\pscircle*[linecolor=lightgray](3.55,1.25){.3}
\pscircle(3.55,1.25){.3}
\pscircle*[linecolor=lightgray](4.65,1.25){.3}
\pscircle(4.65,1.25){.3}
\psline[linewidth=.5pt,linestyle=dashed,dash=3pt
2pt](2.25,1.45)(2.65,1.05)
\psline[linewidth=.5pt,linestyle=dashed,dash=3pt
2pt](2.65,1.47)(2.25,1.05)
\psline[linewidth=.5pt,linestyle=dashed,dash=3pt
2pt](3.35,1.45)(3.75,1.05)
\psline[linewidth=.5pt,linestyle=dashed,dash=3pt
2pt](3.75,1.45)(3.35,1.05)
\psline[linewidth=.5pt,linestyle=dashed,dash=3pt
2pt](4.45,1.45)(4.85,1.05)
\psline[linewidth=.5pt,linestyle=dashed,dash=3pt
2pt](4.85,1.45)(4.45,1.05)
\rput[b](3.6,2){\small$\varepsilon$}
\psline[linewidth=.5pt,linearc=.2](2.65,1.45)(2.75,1.55)(3.25,1.55)(3.35,1.45)
\psline[linewidth=.5pt,linearc=.2](2.65,1.05)(2.75,.95)(3.25,.95)(3.35,1.05)
\psline[linewidth=.5pt,linearc=.2](3.75,1.45)(3.85,1.55)(4.35,1.55)(4.45,1.45)
\psline[linewidth=.5pt,linearc=.2](3.75,1.05)(3.85,.95)(4.35,.95)(4.45,1.05)
\psline[linewidth=.5pt,linearc=.2](2.25,1.45)(2.1,1.6)(2.1,1.9)(3,1.9)
\psline[linewidth=.5pt,linearc=.2](4.85,1.45)(5,1.6)(5,1.9)(3,1.9)
\psline[linewidth=.5pt,linearc=.2](2.25,1.05)(2.1,.9)(2.1,.6)(3,.6)
\psline[linewidth=.5pt,linearc=.2](4.85,1.05)(5,.9)(5,.6)(3,.6)
\psline[linewidth=.5pt, arrowsize=2pt 2.5]{<-}(4.1,2)(4.5,2)
\psline[linewidth=.5pt, arrowsize=2pt 2.5]{->}(4.1,1.8)(4.5,1.8)
\rput[b](1.2,.5){\small$\delta$}
\psframe[linewidth=.5pt,framearc=.2](.7,.4)(3,2.2)
\psline[linewidth=.5pt, arrowsize=2pt 2.5]{->}(.6,1.9)(.6,1.5)
\psline[linewidth=.5pt, arrowsize=2pt 2.5]{->}(.825,1.5)(.825,1.9)
\endpspicture
\end{tabular}
\caption{\label{curves}Generic curves in $F$.}
\end{figure}

Until the end of this paper $F$ will be the non-orientable surface
of genus $4$, obtained by removing from a $2$-sphere four disjoint
open discs and identifying antipodal points on each of the resulting
boundary components. The surface $F$ is represented in Figure
\ref{curves}, where the removed discs are shaded. Let $a_1$, $a_2$,
$a_3$, $a_4$, $b$, $d$ and $e$ denote Dehn twists about the curves
labeled with the corresponding Greek letters in Figure \ref{curves},
in the indicated directions. For $i\in\{1,2,3\}$ we define
$$y_i=y_{\mu_i,\alpha_i},\qquad u_i=a_iy_i.$$
Observe that $u_i$ interchanges $\mu_i$ and $\mu_{i+1}$. We also
define $$t=u_3u_2u_1a_1a_2a_3.$$ A geometric meaning of $t$ will be
explained in Remark \ref{t_geom} below.

We are ready to state our main result:

\begin{theorem}\label{main}
The group $\M(F)$ admits presentation with generators $a_1$, $a_2$,
$a_3$, $a_4$, $b$, $u_1$, $u_2$, $u_3$, $t$ and relations:

\noindent $(1)\ a_1a_3=a_3a_1;\quad$ $(2)\ a_4a_3=a_3a_4;$\\ $(3)$
$ba_1=a_1b$,
$ba_2=a_2b$, $ba_3=a_3b;$\\
$(4)$ $a_1a_2a_1=a_2a_1a_2$, $a_3a_2a_3=a_2a_3a_2$,
$a_4a_2a_4=a_2a_4a_2;$\\
$(5)\ (a_1a_2a_3)^4=1;\quad$ $(6)\ (a_4a_2a_3)^4=1;$\\
$(7)\ u_3a_1u_3^{-1}=a_1;\quad$ $(8)\ u_3a_3u_3^{-1}=a_3^{-1};\quad$
$(9)\
u_3a_2u_3^{-1}=a_2a_4^{-1}a_2^{-1};$\\
$(10)\ (u_3a_4)^2=1;\quad$ $(11)\ (u_3b)^2=1;\quad$ $(12)\
u_3a_4u_3^{-1}=u_1a_4u_1^{-1};$\\ $(13)\ u_1u_3=u_3u_1;\quad$
$(14)\ u_1^2=u_3^2;\quad$ $(15)\ u_1=(a_1a_2a_3)^2u_3(a_1a_2a_3)^2;$\\
$(16)\ u_2=a_3^{-1}a_2^{-1}u_3^{-1}a_2a_3;\quad$ $(17)\
t=u_3u_2u_1a_1a_2a_3;$\\ $(18)\ t^2=1;\quad$ $(19)\
tu_3t=u_3^{-1};\quad$  $(20)\ tbt=b^{-1};\quad$\\ $(21)$
$ta_1=a_1t$, $ta_2=a_2t$, $ta_3=a_3t$.
\end{theorem}

\begin{rem}\label{rem1}
Notice that $a_4$, $u_1$, $u_2$ and $t$ are expressed in terms of
the remaining generators by relations (9,15,16,17). Thus the
presentation in Theorem \ref{main} can be reduced by Tietze
transformations to a presentation with generators $a_1$, $a_2$,
$a_3$, $b$ and $u_3$. This is exactly the generating set for
$\M(F)$ obtained by Chillingworth in \cite{C}. It is not difficult
to show that $\M(F)$ is generated by three elements: $a_1$, $u_3$
and $ba_1a_2a_3$, and it is the minimal size of a generating set
for $\M(F)$ (see \cite{Sz2}).
\end{rem}

%
%



\begin{prop}\label{sat} The relations from Theorem \ref{main}
are satisfied in $\M(F)$.
\end{prop}

\proof Relations (1), (2) and (3) are satisfied, because Dehn twists
about disjoint curves commute. Relations (4) are ``braid'' relations
(\ref{braid}).

Let $\beta'$ and $\beta''$ be boundary curves of a regular
neighborhood of the curve $\beta$, so that $\beta'$ and $\beta''$
also bound a torus with two holes in $F$, which contains the curves
$\alpha_1$, $\alpha_2$ and $\alpha_3$. Then we have ``star''
relation (\ref{star2}): $(a_1a_2a_3)^4=bb^{-1}=1$, hence (5).

Let $\gamma_1$ and $\gamma_2$ be boundary curves of regular
neighborhoods of one-sided curves $\mu_1$, and $\mu_5$, so that
$\gamma_1$ and $\gamma_2$ bound a torus with two holes in $F$, which
contains the curves $\alpha_4$, $\alpha_2$ and $\alpha_3$. From the
``star'' relation (\ref{star2}) we have (6): $(a_4a_2a_3)^4=1$,
because Dehn twists about $\gamma_1$ and $\gamma_2$ are trivial.

Relation (7) is obvious, (8) follows from (\ref{uau}). By (4) we
have $a_2a_4^{-1}a_2^{-1}=a_4^{-1}a_2^{-1}a_4$, hence (9) is
equivalent to $a_4u_3a_2u_3^{-1}a_4^{-1}=a_2^{-1}$ and it can be
verified by checking that $a_4u_3$ fixes $\alpha_2$ and reverses
orientation of its neighborhood.

Let $\gamma_1$ and $\gamma_2$ be boundary curves of regular
neighborhoods of one-sided curves $\mu_1$, and $\mu_2$. Then
$\gamma_1$ and $\gamma_2$ bound a Klein bottle with two holes in $F$
and from (\ref{Krel}) we have (10): $(u_3a_4)^2=1$, because Dehn
twists about $\gamma_1$ and $\gamma_2$ are trivial.

Let $\alpha_1'$ and $\alpha''_1$ be boundary curves of a regular
neighborhood of $\alpha_1$. Then $\alpha_1'$ and $\alpha''_1$ bound
a Klein bottle with two holes in $F$ and from (\ref{Krel}) we have
$(bu_3)^2=a_1a_1^{-1}=1$, hence (11).

Relation (12) can be verified by checking that $u_1^{-1}u_3$ fixes
$\alpha_4$ and preserves orientation of its neighborhood, (13) is
obvious, (14) follows from (\ref{u2}): $u_1^2=d=u_3^2$.

Let $z=(a_1a_2a_3)^{-1}$. It can be checked that
$z(\alpha_3)=\alpha_2^{-1}$, $z(\alpha_2)=\alpha_1^{-1}$ as
oriented curves, and $z(\mu_3)=\mu_2$, $z(\mu_2)=\mu_1$. Hence, by
(\ref{hyh}), we have: $y_2=zy_3^{-1}z^{-1}$ and
$y_1=z^2y_3z^{-2}$. Since $z$ preserves orientation of a regular
neighborhood of $\alpha_1\cup\alpha_2\cup\alpha_3$, thus
$a_2=za_3z^{-1}$ and $a_1=z^2a_3z^{-2}$. Now
$$u_1=a_1y_1=z^2a_3y_3z^{-2}=z^2u_3z^{-2},$$ and since, by (5),
$z^2=z^{-2}=(a_1a_2a_3)^2$, this proves (15). Similarly we prove
(16), using (7) and (8):
$$u_2=a_2y_2=za_3y_3^{-1}z^{-1}=za_3u_3^{-1}a_3z^{-1}\stackrel{\small(8)}{=}zu_3^{-1}z^{-1},$$
$$u_2=a_3^{-1}a_2^{-1}a_1^{-1}u_3^{-1}a_1a_2a_3\stackrel{\small(7)}{=}a_3^{-1}a_2^{-1}u_3^{-1}a_2a_3.$$

Relation (17) is simply definition of $t$. It can be checked, that
for $i\in\{1,2,3\}$, $t$ fixes the curve $\alpha_i$ and preserves
orientation of its neighborhood, hence (21): $ta_it^{-1}=a_i$.
Since $t$ reverses orientation of $\alpha_3$ and fixes $\mu_3$,
thus $ty_3t^{-1}=y_3^{-1}$ and (19):
$$tu_3t^{-1}=ta_3y_3t^{-1}=a_3y_3^{-1}=a_3u_3^{-1}a_3\stackrel{\small(8)}{=}u_3^{-1}.$$
Since $t$ fixes $\beta$ and reverses orientation of its
neighborhood, thus (20): $tbt^{-1}=b^{-1}$. It follows that $t^2$
commutes with $b$, $y_3$ and $a_i$ for $i\in\{1,2,3\}$.  Since these
elements generate $\M(F)$ (see \cite{C}), $t^2$ belongs to the
center of $\M(F)$, which is trivial, according to \cite{S},
Corollary 6.3. Thus (18): $t^2=1$ holds.\foorp

\begin{rem}\label{t_geom}
Recall that $F$ is obtained by removing from a $2$-sphere four
disjoint open discs and identifying antipodal points on each of the
resulting boundary components. Suppose that this sphere is embedded
in $\R^3$, in such a way that it is invariant under reflection about
a plane $\Pi$, which contains centers of the four removed discs
(imagine a plane perpendicular to the plane of Figure \ref{curves},
which contains centers of the four shaded discs). Then, the
reflection about $\Pi$ commutes with the identification, and thus it
induces a homeomorphism of $F$ of order $2$. Denote by $h$ its
isotopy class. It is easy to verify that $ht$ commutes with $b$,
$y_3$ and $a_i$ for $i\in\{1,2,3\}$. Hence we can conclude that
$ht=1$ by arguing as at the end of the proof of Proposition
\ref{sat}. Thus $h=t$. This interpretation of $t$ as being induced
by reflection is convenient for verifying relations involving $t$.
\end{rem}

\medskip

Let $\mathcal{G}$ be an abstract group defined by the presentation
in Theorem \ref{main}. By Proposition \ref{sat}, the map which
assigns to each generator of $\mathcal{G}$ the isotopy class of the
homeomorphism which it represents, extends to a homomorphism
$$\Phi\colon\mathcal{G}\to\M(F).$$
We need to show that $\Phi$ is an isomorphism. Since images of the
generators of $\mathcal{G}$ generate $\M(F)$ (c.f. Remark
\ref{rem1}), $\Phi$ is onto. We will show that it is injective in
Section \ref{inj}.
%

\section{Presentation for $\M(F)$ from its action on $\C$}\label{s_pres}
Recall that the ordered complex of curves $\C$ is a
$\Delta$-complex, whose $n$-simplices are equivalence classes of
generic $(n+1)$-families of disjoint curves. Let $\C^n$ denote the
$n$-skeleton of $\C$, that is the set of its $n$-simplices. Since
generic $n$-families of disjoint curves are ordered $n$-tuples, $\C$
has natural orientation. In particular its edges are oriented. For
an edge $E\in\C^1$ let $i(E)$ and $t(E)$ denote its initial and
terminal vertex respectively. We denote by $\overline{E}$ the
inverse of $E$, that is the edge with the same vertices but opposite
orientation. If $E=[\gamma_1,\gamma_2]$ then $i(E)=[\gamma_1]$,
$t(E)=[\gamma_2]$, $\overline{E}=[\gamma_2,\gamma_1]$.

The mapping class group $\M(F)$ acts on $\C$ by permuting its
simplices,
$h[\gamma_1,\dots,\gamma_n]=[h(\gamma_1),\dots,h(\gamma_n)]$, thus
the orbit space $X=\C/\M(F)$ inherits the structure of a
$\Delta$-complex. Let $X^n$ denote the $n$-skeleton of $X$ and let
$\pi\colon\C\to X$ denote the canonical projection. For $E\in\C^1$
we define $i(\pi(E))=\pi(i(E))$, $t(\pi(E))=\pi(t(E))$,
$\overline{\pi(E)}=\pi(\overline{E})$. We say that $E\in X^1$ is a
{\it loop} based at $V$ if $i(E)=t(E)=V$. In this section we will
define a map $\sigma\colon X^n\to\C^n$ which assigns to each
$n$-simplex of $X$ its representative in $\C$ (i.e.
$\pi\circ\sigma=identity$) for $n=0,1,2$.

\medskip

%
Let $C=(\gamma_1,\dots,\gamma_n)$ be a generic $n$-family of
disjoint curves. Denote by $F_C$ the compact surface obtained by
cutting $F$ along $C$, i.e. the natural compactification of
$F\backslash(\bigcup_{i=1}^{n}\gamma_i)$. Note that $F_C$ is in
general not connected. Denote by $N_1,\dots,N_k$ the connected
components of $F_C$. Then we write
$$\M(F_C)=\M(N_1)\times\dots\times\M(N_k).$$
Denote by $\rho\colon F_C\to F$ the continuous map induced by the
inclusion of $F\backslash(\bigcup_{i=1}^{r}\gamma_i)$ in $F$. The
map $\rho$ induces a homomorphism $\rho_\ast\colon\M(F_C)\to\M(F)$.

Let $\gamma_i$ be a two-sided curve in the family $C$. There exist
two connected components $N'$ and $N''$, and two distinct boundary
curves $\gamma'_i$ and $\gamma''_i$ of $F_C$, such that
$\rho(\gamma'_i)=\rho(\gamma''_i)=\gamma_i$. We say that $\gamma_i$
is a {\it separating limit curve} of $N'$ (and $N''$) if $N'\ne
N''$, and $\gamma_i$ is a {\it non-separating two-sided limit curve}
of $N'$ if $N'=N''$.

Let $\gamma_i$ be a one-sided curve in $C$. There exists a component
$N$  and a boundary curve $\gamma'_i$ of $F_C$ such that
$\rho(\gamma'_i)=\gamma^2_i$. We say that $\gamma_i$ is a {\it
one-sided limit curve} of $N$.

We say that two simplices $[C]$ and $[C']$ of $\C$ are {\it
$\M(F)$-equivalent} if $[C]=h[C']$ for some $h\in\M(F)$. The
following proposition is a special case of Proposition 5.2 of
\cite{Sz1} for closed $F$.

\begin{prop}\label{simplices}
Let $C=(\gamma_1,\dots,\gamma_n)$ and
$C'=(\gamma'_1,\dots,\gamma'_n)$ be two generic $n$-families of
disjoint curves. Then simplices $[C]$ and $[C']$ are
$\M(F)$-equivalent if and only if for all subfamilies $D\subseteq C$
and $D'\subseteq C'$, such that $\gamma_i\in D \iff \gamma'_i\in
D'$, there exists a one to one correspondence between the connected
components of $F_{D}$ and those of $F_{D'}$, such that for every
pair $(N,N')$ where $N$ is any component of $F_{D}$ and $N'$ is the
corresponding component of $F_{D'}$, we have:
\begin{itemize}
\item
$N$ and $N'$ are either both orientable or both non-orientable, of
the same genus;
\item
if $\gamma_i$ is a separating limit curve of $N$, then $\gamma'_i$
is a separating limit curve of $N'$;
\item
if $\gamma_i$ is a non-separating two-sided limit curve of $N$, then
$\gamma'_i$ is a non-separating two-sided limit curve of $N'$;
\item
if $\gamma_i$ is a one-sided limit curve of $N$, then $\gamma'_i$ is
a one-sided limit curve of $N'$.\foorp
\end{itemize}
\end{prop}

\medskip

\begin{prop}\label{X0}
The complex $\Co$ has five $\M(F)$-orbits of vertices represented by
$[\mu_1]$, $[\alpha_3]$, $[\beta]$, $[\delta]$ and $[\varepsilon]$.
\end{prop}
\proof Suppose that $\gamma$ is a non-separating curve in $F$. By
comparing Euler characteristic of $F$ and $F_\gamma$, we obtain that
$F_\gamma$ is non-orientable and has genus $3$ if $\gamma$ is
one-sided, and if $\gamma$ is two-sided, then $F_\gamma$ is either
non-orientable of genus $2$ or orientable of genus $1$. Thus, by
Proposition \ref{simplices}, $\Co$ has three $\M(F)$-orbits of
non-separating vertices, represented by $[\mu_1]$, $[\alpha_3]$ and
$[\beta]$. If $\gamma$ is a separating generic curve, then
$F_\gamma$ is either a disjoint union of two non-orientable surfaces
of genus $2$ or a disjoint union of a non-orientable surface of
genus $2$ and an orientable surface of genus $1$. Thus $\Co$ has two
$\M(F)$-orbits of separating vertices, represented by $[\delta]$ and
$[\varepsilon]$. \foorp

\medskip

By Proposition \ref{X0} the orbit complex $X$ has five vertices. We
denote them by
$$V_1=\pi([\alpha_3]),\ V_2=\pi([\mu_1]),\  V_3=\pi([\beta]),\
V_4=\pi([\delta]),\ V_5=\pi([\varepsilon]).$$ We also define a
section $\sigma\colon X^0\to\C^0$ by
$$\sigma(V_1)=[\mu_1],\  \sigma(V_2)=[\alpha_1],\  \sigma(V_3)=[\beta],\
\sigma(V_4)=[\delta],\  \sigma(V_5)=[\varepsilon].$$ For each
$V\in X^0$ let $S_V=\Stab(\sigma(V))$ denote the stabilizer of
$\sigma(V)$ in $\M(F)$.

\medskip


\begin{table}
\caption{\label{ch_edges}Edges.}
\begin{tabular}{|c|c|c|c|c|}
\hline
$E$&$\sigma(E)$&$\sigma(t(E))$&$g_E$&$G_E$\\
\hline $E_1$&$[\alpha_3,\mu_1]$&$[\mu_1]$&$1$&$\{a_3,a_4,u_3,t,y_1\}$\\
\hline
$E_2$&$[\alpha_3,\beta]$&$[\beta]$&$1$&$\{b,a_1,a_3,(a_3^2a_2)^2,t,u_1^{-1}u_3\}$\\
\hline $E_3$&$[\alpha_3,\delta]$&$[\delta]$&$1$&$\{a_3,a_1,u_1,u_3,t\}$\\
\hline
$E_4$&$[\alpha_3,\varepsilon]$&$[\varepsilon]$&$1$&\\
\hline
$E_5$&$[\beta,\varepsilon]$&$[\varepsilon]$&$1$&\\
\hline $E_6$&$[\mu_1,\varepsilon]$&$[\varepsilon]$&$1$&$\{a_2,a_3,t,u_3u_2u_3\}$\\
\hline $E_7$&$[\mu_1,\delta]$&$[\delta]$&$1$&$\{a_3,u_3,t,y_1\}$\\
\hline $E_8$&$[\alpha_3,\alpha_1]$&$[\alpha_3]$&$(a_1a_2a_3)^2$&$\{a_1,a_3,b,u_1,u_3,t\}$\\
\hline
$E_9$&$[\alpha_3,\alpha_4]$&$[\alpha_3]$&$a_2^{-1}u_2^{-1}$&$\{a_3,a_4,u_3b,u_1b,u_1t\}$\\
\hline $E_{10}$&$[\mu_1,\mu_2]$&$[\mu_1]$&$u_1$&$\{u_3,a_3,a_4,t,y_2\}$\\
\hline $E_{11}$&$[\mu_1,\mu_5]$&$[\mu_1]$&$b^{-1}$&$\{a_2,a_3,a_4,u_3u_2u_3t\}$\\
\hline
\end{tabular}
\end{table}

For $i\in\{1,\dots,11\}$ we define an edge $E_i\in X^1$ by
$E_i=\pi(\sigma(E_i))$, where $\sigma(E_i)$ is an edge of $\Co$
defined in the second column of Table \ref{ch_edges}.

\begin{prop}\label{X1} Every edge of $\Co$ is $\M(F)$-equivalent to
$\sigma(E_i)$ or $\overline{\sigma(E_i)}$ for some
$i\in\{1,\dots,11\}$.
\end{prop}
\proof Let $(\gamma_1,\gamma_2)$ be a generic pair of disjoint
curves representing an edge of $\C$. By Proposition \ref{X0},
$[\gamma_i]$ is $\M(F)$-equivalent to one of the vertices $[\mu_1]$,
$[\alpha_1]$, $[\beta]$, $[\delta]$ or $[\varepsilon]$.

Suppose that $[\gamma_2]$ is $\M(F)$-equivalent to $[\delta]$. Then
$F_{\gamma_2}$ has two connected components, each homeomorphic to
the Klein bottle with a hole. Denote by $N$ the component containing
$\gamma_1$. If $\gamma_1$ is one-sided, then $N_{\gamma_1}$ is
projective plane with two holes. If $\gamma_1$ is two-sided, then
since it is generic and not isotopic to $\gamma_2$, it is
non-separating, $N_{\gamma_1}$ is pair of pants and $F_{\gamma_1}$
is non-orientable. Thus by Proposition \ref{simplices},
$[\gamma_1,\gamma_2]$ is $\M(F)$-equivalent to $\sigma(E_3)$ or
$\sigma(E_7)$.

Suppose that $[\gamma_2]$ is $\M(F)$-equivalent to $[\varepsilon]$.
Then $F_{\gamma_2}$ has components $N$ homeomorphic to the Klein
bottle with a hole and $N'$ homeomorphic to the torus with a hole.
If $\gamma_1\subset N$, then as above, $N_{\gamma_1}$ is projective
plane with with two holes if $\gamma_1$ is one-sided, or pair of
pants if it is two-sided. If $\gamma_1$ is two-sided then
$F_{\gamma_1}$ is orientable. If $\gamma_1\subset N'$, then
$\gamma_1$ is two-sided and non-separating, $N'_{\gamma_1}$ is pair
of pants and $F_{\gamma_1}$ is non-orientable. Thus by Proposition
\ref{simplices}, $[\gamma_1,\gamma_2]$ is $\M(F)$-equivalent to
$\sigma(E_4)$ or $\sigma(E_5)$ or $\sigma(E_6)$.

If $\gamma_1$ is separating, then clearly $[\gamma_1,\gamma_2]$ is
$\M(F)$-equivalent to $\overline{\sigma(E_i)}$ for some
$i\in\{3,\dots,7\}$. It remains to consider cases where $\gamma_i$
are non-separating.

Suppose that $[\gamma_2]$ is $\M(F)$-equivalent to $[\beta]$. Then
$F_{\gamma_2}$ is torus with two holes. Since $\gamma_1$ is
non-separating in $F$ and not isotopic to $\gamma_2$, thus it is
also non-separating in $F_{\gamma_2}$ and $F_{(\gamma_1,\gamma_2)}$
is sphere with four holes. Note that $F_{\gamma_1}$ is
non-orientable, thus by Proposition \ref{simplices},
$[\gamma_1,\gamma_2]$ is $\M(F)$-equivalent to $\sigma(E_2)$

Suppose that $[\gamma_2]$ is $\M(F)$-equivalent to $[\alpha_3]$.
Then $F_{\gamma_2}$ is Klein bottle with two holes. If $\gamma_1$ is
one-sided, then $F_{(\gamma_1,\gamma_2)}$ is projective plane with
$3$ holes and $[\gamma_1,\gamma_2]$ is $\M(F)$-equivalent to
$\overline{\sigma(E_1)}$. Suppose that $\gamma_1$ is two-sided. If
it is non-separating in $F_{\gamma_2}$, then
$F_{(\gamma_1,\gamma_2)}$ is sphere with $4$ holes and
$[\gamma_1,\gamma_2]$ is $\M(F)$-equivalent to $\sigma(E_8)$ if
$F_{\gamma_1}$ is non-orientable, or to $\overline{\sigma(E_2)}$ if
$F_{\gamma_1}$ is orientable. If  $\gamma_1$ is separating in
$F_{\gamma_2}$ (but non-separating in $F$), then
$F_{(\gamma_1,\gamma_2)}$ is disjoint union of two copies of the
projective plane with two holes and $F_{\gamma_1}$ is
non-orientable. Thus $[\gamma_1,\gamma_2]$ is $\M(F)$-equivalent to
$\sigma(E_9)$.

It remains to consider the case when $\gamma_i$ are one-sided. Then
$F_{(\gamma_1,\gamma_2)}$ is connected and if it is non-orientable,
then $[\gamma_1,\gamma_2]$ is $\M(F)$-equivalent to
$\sigma(E_{10})$. Otherwise $[\gamma_1,\gamma_2]$ is
$\M(F)$-equivalent to $\sigma(E_{11})$. \foorp

\medskip

Since for $8\le j\le 11$ the edges $\sigma(E_j)$ and
$\overline{\sigma(E_j)}$ are $\M(F)$-equivalent, hence
$\overline{E_j}=E_j$. Thus Proposition \ref{X1} asserts that
$$X^1=\{E_i, \overline{E_j}\,|\,1\le i\le 11, 1\le j\le 7\}.$$
For $i\in\{1,\dots,7\}$ we define
$\sigma(\overline{E_i})=\overline{\sigma(E_i)}$. For each $E\in X^1$
let $S_E=\Stab(\sigma(E))$.

Observe that for each $E\in X^1$ we have
$i(\sigma(E))=\sigma(i(E))$. For $i\in\{1,\dots,11\}$ let $g_{E_i}$
be the element of $\M(F)$ defined in the fourth column of Table
\ref{ch_edges}. For $j\in\{1,\dots,7\}$ let $g_{\overline{E_j}}=1$.
It can be checked that for each $E\in X^1$
$$g_E(\sigma(t(E)))=t(\sigma(E)).$$ The conjugation map $c_E$
defined by $g\mapsto g_E^{-1}gg_E$ maps $\Stab(t(\sigma(E)))$ onto
$\Stab(\sigma(t(E)))$; in particular, $c_E(S_E)\subseteq S_{t(E)}$.

We define
$$\T=\{E_1, E_2, E_3, E_4\}.$$
Note that $\mathcal{T}$ is a maximal tree in $X^1$ regarded as a
graph.

\medskip


\begin{table}
\caption{\label{traingles}Triangles.}
\begin{tabular}{|c|c|c|}
\hline
$T$&$\sigma(T)$& edges\\
\hline $T_1$&$[\alpha_3,\mu_1,\mu_2]$&$E_1$, $E_{10}$, $E_1$\\
\hline
$T_2$&$[\alpha_3,\mu_1,\mu_5]$&$E_1$, $E_{11}$, $E_1$\\
\hline $T_3$&$[\alpha_3,\mu_1,\delta]$&$E_1$, $E_7$, $E_3$\\
\hline
$T_4$&$[\alpha_3,\alpha_4,\mu_1]$&$E_9$, $E_1$, $E_1$\\
\hline
$T_5$&$[\alpha_3,\mu_1,\varepsilon]$&$E_1$, $E_6$, $E_4$\\
\hline $T_6$&$[\alpha_3,\alpha_1,\beta]$&$E_8$, $E_2$, $E_2$\\
\hline $T_7$&$[\alpha_3,\beta,\varepsilon]$&$E_3$, $E_5$, $E_4$\\
\hline $T_8$&$[\alpha_3,\alpha_1,\delta]$&$E_8$, $E_3$, $E_3$\\
\hline
$T_9$&$[\mu_1,\mu_5,\varepsilon]$&$E_{11}$, $E_6$, $E_6$\\
\hline $T_{10}$&$[\mu_1,\mu_3,\delta]$&$E_{10}$, $E_7$, $E_7$\\
\hline $T_{11}$&$[\mu_1,\mu_2,\delta]$&$E_{10}$, $E_7$, $E_7$\\
\hline $T_{12}$&$[\mu_1,\mu_2,\mu_3]$&$E_{10}$, $E_{10}$, $E_{10}$\\
\hline
\end{tabular}
\end{table}

For $i\in\{1,\dots,12\}$ we define a triangle $T_i\in X^2$ by
$T_i=\pi(\sigma(T_i))$, where $\sigma(T_i)$ is a triangle of $\Co$
defined in the second column of Table \ref{traingles}.

\begin{prop}\label{X2}
Let $(\gamma_1,\gamma_2,\gamma_3)$ be any generic $3$-family of
disjoint curves in $F$. Then there exists a permutation
$\tau\in\Sigma_3$ such that the simplex
$[\gamma_{\tau(1)},\gamma_{\tau(2)},\gamma_{\tau(3)}]$ of $\Co$ is
$\M(F)$-equivalent to $\sigma(T_i)$ for some $i\in\{1,\dots,12\}$.
\end{prop}
\proof Let $T=\pi([\gamma_1,\gamma_2,\gamma_3])$,
$A=\pi([\gamma_1,\gamma_2])$, $B=\pi([\gamma_2,\gamma_3])$,
$C=\pi([\gamma_1,\gamma_3])$.

Suppose that at least one edge of $T$ is $E_1$. By permuting the
vertices of $T$ we may assume that $A=E_1$. Then
$[\gamma_1,\gamma_2]$ is $\M(F)$-equivalent to
$\sigma(E_1)=[\alpha_3,\mu_1]$ and $F_{(\gamma_1,\gamma_2)}$ is
projective plane with $3$ holes.

Suppose that $\gamma_3$ is one-sided. Then
$F_{(\gamma_1,\gamma_2,\gamma_3)}$ is sphere with four holes and
$C=E_1$. If $F_{(\gamma_2,\gamma_3)}$ is non-orientable, then by
Proposition \ref{simplices}, $B=E_{10}$ and $T=T_1$. Otherwise
$B=E_{11}$ and $T=T_2$.

Suppose that $\gamma_3$ is separating. Then
$F_{(\gamma_1,\gamma_2,\gamma_3)}$ is disjoint union of par of pants
and projective plane with two holes. If both components of
$F_{\gamma_3}$ are non-orientable, that is if $[\gamma_3]$ is
$\M(F)$-equivalent to $[\delta]$, then $B=E_7$, $C=E_3$ and $T=T_3$.
If one component of $F_{\gamma_3}$ is orientable, that is if
$[\gamma_3]$ is $\M(F)$-equivalent to $[\varepsilon]$, then $B=E_6$,
$C=E_4$ and $T=T_5$.

Suppose that $\gamma_3$ is two-sided and non-separating, that is
$[\gamma_3]$ is $\M(F)$-equivalent to $[\alpha_3]$. Then it must be
separating in $F_{(\gamma_1,\gamma_2)}$ and
$F_{(\gamma_1,\gamma_2,\gamma_3)}$ is again disjoint union of par of
pants and projective plane with two holes. By Proposition
\ref{simplices}, $B=\overline{E_1}$, $C=E_9$ and
$\pi([\gamma_1,\gamma_3,\gamma_2])=T_4$.

Suppose that at least one edge if $T$ is $E_2$. By permuting the
vertices of $T$ we may assume that $A=E_2$. Then
$[\gamma_1,\gamma_2]$ is $\M(F)$ equivalent to
$\sigma(E_2)=[\alpha_3,\beta]$ and $F_{(\gamma_1,\gamma_2)}$ is
sphere with $4$ holes. Now $\gamma_3$ is two-sided and
$F_{(\gamma_1,\gamma_2,\gamma_3)}$ is disjoint union of two pairs of
pants. If $\gamma_3$ is separating in $F$, then $[\gamma_3]$ is
$\M(F)$-equivalent to $[\varepsilon]$, $B=E_5$, $C=E_4$ and $T=T_7$.
If $\gamma_3$ is non-separating, then $[\gamma_3]$ is
$\M(F)$-equivalent to $[\alpha_3]$, $B=\overline{E_2}$, $C=E_8$ and
$\pi([\gamma_1,\gamma_3,\gamma_2])=T_6$.

For the rest of the proof we may assume that no edge of $T$ is equal
to $E_1$, $E_2$, $\overline{E_1}$ or $\overline{E_2}$. Suppose that
$[\gamma_1]$ is $\M(F)$-equivalent to $[\alpha_3]$. Since there is
no edge in $\Co$ between two separating vertices, $\gamma_2$ or
$\gamma_3$ must be non-separating. By permuting the vertices we may
assume that is $\gamma_2$. Thus $A=E_8$ or $A=E_9$. Suppose $A=E_8$.
Then $F_{(\gamma_1,\gamma_2)}$ is sphere with $4$ holes and
$F_{(\gamma_1,\gamma_2,\gamma_3)}$ is disjoint union of two pairs of
pants. Note that $\gamma_3$ must be separating, because otherwise
$[\gamma_3]$ would be $\M(F)$ equivalent to $[\beta]$ and $C=E_2$,
which contradicts our assumption about the edges of $T$. Thus
$[\gamma_3]$ is $\M(F)$ equivalent to $[\delta]$, $B=C=E_3$ and
$T=T_8$. Suppose $A=E_9$. Then $F_{(\gamma_1,\gamma_2)}$ is disjoint
union of two copies of projective plane with two holes. But then
$\gamma_3$ must be one-sided and $C=E_1$, which also contradicts the
assumption about the edges of $T$.

For the rest of the proof we assume that no vertex of $T$ is equal
to $\pi[\alpha_3]$. Since there is no edge between two separating
vertices and there is no loop at $\pi([\beta])$, at least one vertex
of $T$ is one-sided. But there is no edge between $\pi([\beta])$ and
$\pi([\mu_1])$, hence no vertex of $T$ is equal to $\pi[\beta]$.
Thus $T$ has at least two one-sided vertices. By permuting the
vertices of $T$ we may assume that $\gamma_1$ and $\gamma_2$ are
one-sided, hence $A\in\{E_{10},E_{11}\}$.

Suppose that $A=E_{10}$. Then $F_{(\gamma_1,\gamma_2)}$ is Klein
bottle with two holes. If $\gamma_3$ is one-sided, then
$F_{(\gamma_1,\gamma_2,\gamma_3)}$ is non-orientable and $T=T_{12}$.
If $\gamma_3$ is separating, then it is $\M(F)$ equivalent to
$[\delta]$. If $\gamma_1$ and $\gamma_2$ are in the same component
of $F_{\gamma_3}$ then $T=T_{11}$. Otherwise $T=T_{10}$.

Suppose that $A=E_{11}$. Then $F_{(\gamma_1,\gamma_2)}$ is torus
with two holes, $\gamma_3$ is separating $\M(F)$-equivalent to
$[\varepsilon]$ and $T=T_9$.
 \foorp

\medskip

Proposition \ref{X2} asserts that
$$X^2=\{T^\tau_i\,|\,i\in\{1,\dots,12\}, \tau\in\Sigma_3\},$$
where
$T^\tau_i=\pi([\gamma_{\tau(1)},\gamma_{\tau(2)},\gamma_{\tau(3)}])$
if $T_i=\pi([\gamma_1,\gamma_2,\gamma_3])$. Observe that
$T^\tau_{12}=T_{12}$ for each $\tau\in\Sigma_3$, for $i\in\{3,5,7\}$
permutations of vertices yield $6$ different triangles $T^\tau_i$,
whereas for $i\ne 3,5,7,12$ there are $3$ different triangles
$T^\tau_i$. For example for $i=1$ these are:
$$T^1_1=\pi([\alpha_3,\mu_1,\mu_2]),\ T^{(1,2)}_1=\pi([\mu_1,\alpha_3,\mu_2]),\ T^{(1,3)}_1=\pi([\mu_1,\mu_2,\alpha_3]).$$

\begin{figure}
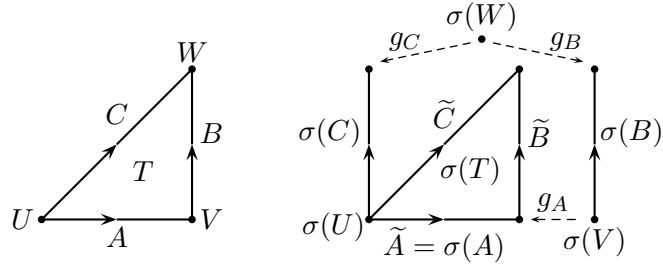

\begin{tabular}{cc}
\input{tri1}&\input{tri2}
\end{tabular}
\caption{\label{tr}Triangle in $X$ and its representative in $\C$.}
\end{figure}

For every triangle $T=T^\tau_i\in X^2$, with edges $A$, $B$, $C$
such that $i(C)=i(A)=U$, $t(A)=i(B)=V$, $t(B)=t(C)=W$, we choose a
representative $\sigma(T)$ in $\C^2$ by permuting vertices of
$\sigma(T_i)$. Notice that we can always do it in such a way that if
$\widetilde{A}$, $\widetilde{B}$, $\widetilde{C}$ are the
corresponding edges of $\sigma(T)$, then
$i(\widetilde{C})=i(\widetilde{A})=\sigma(U)$ and
$\widetilde{A}=\sigma(A)$ (see Figure \ref{tr}). For example for
$i=1$:
$$\sigma(T^1_1)=\sigma(T_1),\ \sigma(T^{(1,2)}_1)=[\mu_1,\alpha_3,\mu_2],\
\sigma(T^{(1,3)}_1)=[\mu_1,\mu_2,\alpha_3].$$

Then we can choose elements
$$\varphi\in S_V,\quad \psi\in S_W,\quad \eta\in S_U,$$
such that  $g_A\varphi(\sigma(B))=\tilde{B}$, $g_A\varphi g_B\psi
g^{-1}_C(\sigma(C))=\tilde{C}$, $\eta=g_A\varphi g_B\psi g^{-1}_C.$


\medskip

The next theorem is a special case of a general result of Brown
\cite{Br} (c.f. Theorem 6.3 of \cite{Sz1}).

\begin{theorem}\label{presentation}
Suppose that:

(1) for each $V\in X^0$ the stabilizer $S_V$ has presentation
$\langle G_V|R_V\rangle$;

(2) for each $E\in X^1$ the stabilizer $S_E$ is generated by $G_E$;

Then $\M(F)$ admits the presentation:
\begin{eqnarray*}
\textrm{generators}&=&\bigcup_{V\in X^0}G_V\cup
\{g_E\,|\,E\in X^1\},\\
\textrm{relations}&=& \bigcup_{V\in X^0}R_V \cup R^{(1)}\cup
R^{(2)}\cup R^{(3)},
\end{eqnarray*}
where:\\
$R^{(1)}:\quad g_E=1$ for $E\in\mathcal{T}$;\\
$R^{(2)}:\quad g_E^{-1}i_E(g)g_E=c_E(g)$ for $E\in X^1, g\in G_E$,
where $i_E$ is the inclusion $S_E\hookrightarrow S_{i(E)}$ and
$c_E\colon S_E\to S_{t(E)}$ is the conjugation map defined above;\\
$R^{(3)}:\quad g_A\varphi g_B\psi g^{-1}_C=\eta$ for $T\in X^2$.

\medskip

\noindent In $R^{(2)}$ and $R^{(3)}$, $i_E(g)$, $c_E(g)$, $\varphi$,
$\psi$ and $\eta$ should be expressed as words in the generators
$\bigcup_{V\in X^0}G_V.$ \foorp\end{theorem}

\section{\label{s_stab}Stabilizers of vertices and edges}

Let $C=(\gamma_1,\dots,\gamma_n)$ be a generic $n$-family of
disjoint curves. The stabilizer $\textrm{Stab}[C]$ consist of the
isotopy classes of all homeomorphisms fixing each curve of $C$ (see
\cite{Sz1}). Let $\textrm{Stab}^+[C]$ denote its subgroup consisting
of the isotopy classes of homeomorphism which also preserve the
orientation of each curve of $C$. Clearly $\textrm{Stab}^+[C]$ is a
normal subgroup of $\textrm{Stab}[C]$ of index at most $2^n$.
Observe that the image of $\rho_\ast\colon\M(F_C)\to\M(F)$ is
contained in $\textrm{Stab}^+[C]$ and it consists of the isotopy
classes of homeomorphisms which preserve orientation of a regular
neighborhood of each two-sided curve of $C$ (equivalently they
preserve sides of such curve).

For each curve $\gamma_i\in C$ we define an element
$k_i\in\ker\rho_\ast$ as follows. If $\gamma_i$ is one-sided, then
let $\gamma'_i$ denote the boundary curve of $F_C$, such that
$\rho(\gamma'_i)=\gamma^2_i$. We define $k_i$ to be a Dehn twist
about $\gamma'_i$. If $\gamma_i$ is two-sided, then let $\gamma'_i$
and $\gamma''_i$ denote the boundary curves of $F_C$, such that
$\rho(\gamma'_i)=\rho(\gamma''_i)=\gamma_i$. Let $c'_i$ and $c''_i$
be Dehn twists about these boundary curves, such that
$\rho_\ast(c'_i)=\rho_\ast(c''_i)$. Then we define
$k_i=c'_i(c''_i)^{-1}$. The subgroup of $\M(F_C)$ generated by
$k_1,\dots,k_n$ is a free abelian group of rank $n$ (by \cite{S},
Proposition 4.4) and is equal to $\ker\rho_\ast$ by \cite{Sz1},
Lemma 4.1. Hence we have an exact sequence
\begin{equation}\label{sequence}
1\to\Z^n\to\M(F_C)\stackrel{\rho_\ast}{\to}\textrm{Stab}^+[C]\to\Z_2^r,
\end{equation}
where $r$ is the number of two-sided curves in $C$. By using
sequence (\ref{sequence}) wy may determine a presentation for
$\Stab^+[C]$, and then also for $\Stab[C]$, starting from a
presentation for $\M(F_C)$.

\begin{prop}\label{stab_mu1}
The stabilizer $S_{V_2}=\mathrm{Stab}[\mu_1]$ admits a presentation
with generators $a_2$, $a_3$, $a_4$, $u_2$, $u_3$, $t$ and
relations:

\noindent $\mathrm{(i)}\ a_3a_4=a_4a_3,\quad$ $\mathrm{(ii)}\
a_2a_3a_2=a_3a_2a_3,\quad$ $\mathrm{(iii)}\
a_2a_4a_2=a_4a_2a_4,\quad$\\ $\mathrm{(iv)}\
u_3a_3u_3^{-1}=a_3^{-1},\quad$ $\mathrm{(v)}\
u_3a_2u_3^{-1}=a_2a_4^{-1}a_2^{-1},\quad$ $\mathrm{(vi)}\
(u_3a_4)^2=1,\quad$\\ $\mathrm{(vii)}\ (a_4a_2a_3)^4=1,\quad$
$\mathrm{(viii)}\ t^2=1,\quad$ $\mathrm{(ix)}\
tu_3t=u_3^{-1},\quad$\\ $\mathrm{(x)}\ ta_2=a_2t,\quad$
$\mathrm{(xi)}\ ta_3=a_3t,\quad$ $\mathrm{(xii)}\
u_2=a_3^{-1}a_2^{-1}u_3^{-1}a_2a_3,\quad$\\ $\mathrm{(xiii)}\
u_2a_2u_2^{-1}=a_2^{-1},\quad$ $\mathrm{(xiv)}\ tu_2t=u_2^{-1}.$

\noindent Relations \textrm{(i--xiv)} are consequences of relations
\textrm{(1--21)} in Theorem \ref{main}.
\end{prop}

\proof Notice that (i-xii) appear among relations (1--21) in Theorem
\ref{main}. We will show that $\mathrm{Stab}[\mu_1]$ admits a
presentation with generators $a_2$, $a_3$, $a_4$, $u_3$, $t$ and
relations (i--xi). By Theorem 7.16 of \cite{Sz1} the group
$\M(F_{\mu_1})$ admits a presentation with generators $a_2$, $a_3$,
$a_4$, $u_3$ and relations (i--v) and
$$(u_3a_4)^2=(a_4u_3)^2=(a_4a_2a_3)^4.$$ By (\ref{Krel}), $(u_3a_4)^2$ is
a Dehn twist about $\bdr{F_{\mu_1}}$, hence it
generates the kernel of
$\rho_\ast\colon\M(F_{\mu_1})\to\mathrm{Stab}^+[\mu_1]$. Since
$\rho_\ast$ is onto,
$$\mathrm{Stab}^+[\mu_1]=\langle a_2, a_3, a_4, u_3\,|\,\mathrm{(i-vii)}\rangle.$$
Observe that $t$ reverses the orientation of $\mu_1$ and hence it
represents the nontrivial coset of $\mathrm{Stab}^+[\mu_1]$ in
$\mathrm{Stab}[\mu_1]$. It follows that the last group is generated
by $a_2$, $a_3$, $a_4$, $u_3$ and $t$ satisfying as defining
relations (i--vii), $t^2\in\mathrm{Stab}^+[\mu_1]$ and
$tht\in\mathrm{Stab}^+[\mu_1]$, for $h\in\{a_2,a_3,a_4,u_3\}$.
Notice that (viii--xi) have this form and they hold in $\M(F)$ by
Proposition \ref{sat}. Finally notice that
$ta_4t\in\mathrm{Stab}^+[\mu_1]$ is a consequence of (v)
$a_4=a_2^{-1}u_3a_2^{-1}u_3^{-1}a_2$ and (ix),(x).

Now it remains to show that relations (xiii) and (iv) hold in
$\M(F)$. Indeed, (xiii) follows from (\ref{uau}), while (xiv) is an
easy consequence of (16,18,19,21) in Theorem \ref{main}. Since
(i-xi) are defining relations for $\mathrm{Stab}[\mu_1]$, (xiii) is
a consequence of (i-xii), hence also of (1--21).
\foorp

\begin{prop}\label{stab_d}
The stabilizer $S_{V_4}=\mathrm{Stab}[\delta]$ admits a presentation
with generators $a_1$, $a_3$, $u_1$, $u_3$, $s=(a_1a_2a_3)^2$, and
relations:

\noindent $\mathrm{(i)}\ u_1a_1u_1^{-1}=a_1^{-1},\quad$
$\mathrm{(ii)}\ u_3a_3u_3^{-1}=a_3^{-1},\quad$ $\mathrm{(iii)}\
u_1^2=u_3^2,\quad$\\ $\mathrm{(iv)}\ u_1u_3=u_3u_1,\quad$
$\mathrm{(v)}\ a_1u_3=u_3a_1,\quad$ $\mathrm{(vi)}\
u_1a_3=a_3u_1,\quad$\\ $\mathrm{(vii)}\ a_1a_3=a_3a_1,\quad$
$\mathrm{(viii)}\ t^2=1,\quad$ $\mathrm{(ix)}\ ta_1=a_1t,\quad$
$\mathrm{(x)}\ ta_3=a_3t,\quad$\\ $\mathrm{(xi)}\
tu_1t=u_1^{-1},\quad$ $\mathrm{(xii)}\ tu_3t=u_3^{-1},\quad$
$\mathrm{(xiii)}\ s^2=1,\quad$\\ $\mathrm{(xiv)}\ sa_1s=a_3,\quad$
$\mathrm{(xv)}\ su_1s=u_3,\quad$ $\mathrm{(xvi)}\ st=ts.$

\noindent Relations \textrm{(i--xvi)} are consequences of relations
\textrm{(1--21)} in Theorem \ref{main}.
\end{prop}

\proof First we show that (i--xvi) are consequences of (1--21).
Notice that (i), (vi) and (xi) follow easily from (ii), (v) and
(xii--xv). Relations (ii,iii,v,vii--x,xii) appear among relations
(1--21) in Theorem \ref{main}; (xiii) and (xv) are (5) and (15)
respectively. Relations (1) and (4) imply $sa_1=a_3s$, which
together with (xiii) gives (xiv). Finally, (xvi) follows from (21).

The surface $F_\delta$ has two connected components, each
homeomorphic to the Klein bottle with a hole. By Theorem A.7 of
\cite{S} we have
$$\M(F_\delta)=\langle a_1,
u_1\,|\,u_1a_1u_1^{-1}=a_1^{-1}\rangle\times\langle a_3,
u_3\,|\,u_3a_3u_3^{-1}=a_3^{-1}\rangle.$$ By (\ref{u2}),
$u_1^2=u_3^2=d$, hence $\ker\rho_\ast$ is generated by
$u_1^2u_3^{-2}$ and
$$\rho_\ast(\M(F_\delta))=\langle a_1,a_3,u_1,u_3\,|\,\mathrm{(i-vii)}\rangle.$$
Observe that $s$ and $t$ fix $\delta$ and reverse its orientation,
$s$ preserves, while $t$ reverses orientation of its regular
neighborhood. It follows that (i-xvi) are defining relations for
$\mathrm{Stab}[\delta]$. \foorp

\begin{prop}\label{stab_a3}
The stabilizer $S_{V_1}=\mathrm{Stab}[\alpha_3]$ admits a
presentation with generators $a_1$, $a_3$, $a_4$, $b$ $u_1$, $u_3$,
$t$ and relations:

\noindent $\mathrm{(i)}\ a_1b=ba_1,\quad$ $\mathrm{(ii)}\
u_1a_1u_1^{-1}=a_1^{-1},\quad$ $\mathrm{(iii)}\
ba_4b^{-1}=u_1^{-1}a_4^{-1}u_1,\quad$\\ $\mathrm{(iv)}\
(u_1b)^2=1,\quad$ $\mathrm{(v)}\ (u_1a_4)^2=1,\quad$ $\mathrm{(vi)}\
a_3b=ba_3,\quad$\\ $\mathrm{(vii)}\ a_1a_3=a_3a_1,\quad$
$\mathrm{(viii)}\ a_3a_4=a_4a_3,\quad$ $\mathrm{(ix)}\
a_3u_1=u_1a_3,\quad$\\ $\mathrm{(x)}\ u_3^2=u_1^2,\quad$
$\mathrm{(xi)}\ u_3a_1=a_1u_3,\quad$ $\mathrm{(xii)}\
u_3a_3u_3^{-1}=a_3^{-1},\quad$\\ $\mathrm{(xiii)}\
u_3bu_3^{-1}=u_1bu_1^{-1},\quad$ $\mathrm{(xiv)}\
u_3a_4u_3^{-1}=u_1a_4u_1^{-1},\quad$\\ $\mathrm{(xv)}\
u_3u_1=u_1u_3,\quad$ $\mathrm{(xvi)}\ t^2=1,\quad$ $\mathrm{(xvii)}\
ta_1=a_1t,\quad$\\ $\mathrm{(xviii)}\ ta_3=a_3t,\quad$
$\mathrm{(xix)}\ ta_4t=u_1^{-1}a_4^{-1}u_1,\quad$ $\mathrm{(xx)}\
tbt=b^{-1},\quad$\\ $\mathrm{(xxi)}\ tu_1t=u_1^{-1},\quad$
$\mathrm{(xxii)}\ tu_3t=u_3^{-1}$.

\noindent Relations \textrm{(i--xxii)} are consequences of relations
\textrm{(1--21)} in Theorem \ref{main}.
\end{prop}

\proof First we show that (i--xxii) are consequences of (1--21).
Relations (i,vi--viii,x--xii,xiv-xviii,xx,xxii) appear among
relations (1--21) in Theorem \ref{main}, while (ii,ix,xxi) appear in
Proposition \ref{stab_d}. Relation (iv) follows from (3,5,11,15):
$$(u_1b)^2\stackrel{\small(5,15)}{=}((a_1a_2a_3)^{-2}u_3(a_1a_2a_3)^2b)^2\stackrel{(3)}{=}(a_1a_2a_3)^{-2}(u_3b)^2(a_1a_2a_3)^2
\stackrel{(11)}{=}1.$$ Relation (v) follows from (10,12,14):
$$(u_1a_4)^2=u_1a_4u_1^{-1}u_1^2a_4\stackrel{(12,14)}{=}u_3a_4u_3^{-1}u_3^2a_4=(u_3a_4)^2\stackrel{(10)}{=}1.$$
 Relation (xiii) follows from (11,14) and (iv):
 $$u_3bu_3^{-1}\stackrel{\small(11)}{=}b^{-1}u_3^{-2}\stackrel{(14)}{=}b^{-1}u_1^{-2}\stackrel{\small\mathrm{(iv)}}{=}
 u_1bu_1^{-1}.$$
 By (9) we have $a_4=a_2^{-1}u_3a_2^{-1}u_3^{-1}a_2$, and by (3,11)
 $$ba_4b^{-1}=ba_2^{-1}u_3a_2^{-1}u_3^{-1}a_2b^{-1}=a_2^{-1}bu_3a_2^{-1}u_3^{-1}b^{-1}a_2
 =a_2^{-1}u_3^{-1}a_2^{-1}u_3a_2.$$
 Since, by (12), $u_1^{-1}a_4^{-1}u_1=u_3^{-1}a_4^{-1}u_3$, (iii) is
 equivalent to
$$a_2^{-1}u_3^{-1}a_2^{-1}u_3a_2=u_3^{-1}a_4^{-1}u_3\Leftrightarrow
u_3a_2^{-1}u_3^{-1}a_2^{-1}u_3a_2u_3^{-1}a_4=1.$$ The last relation
is a consequence of (4,9):
$$u_3a_2^{-1}u_3^{-1}a_2^{-1}u_3a_2u_3^{-1}a_4\stackrel{\small(9)}{=}
a_2a_4a_2^{-1}a_4^{-1}a_2^{-1}a_4\stackrel{\small(4)}{=}1.$$
Finally, from (18,19,21) we have:
$$ta_4t=ta_2^{-1}u_3a_2^{-1}u_3^{-1}a_2t=a_2^{-1}u_3^{-1}a_2^{-1}u_3a_2=ba_4b^{-1}\stackrel{\small\mathrm{(iii)}}{=}u_1a_4^{-1}u_1^{-1},$$
that is relation (xix).

The surface $F_{\alpha_3}$ is Klein bottle with two holes. Let
$a'_3$, $a''_3$ denote Dehn twists about its boundary components,
such that $\rho_\ast(a'_3)=\rho_\ast(a''_3)=a_3$. Then, by Theorem
7.10 of \cite{Sz1}, $\M(F_{\alpha_3})$ admits a presentation with
generators $a_1$, $a_4$, $b$ $u_1$, $a'_3$, $a''_3$ and relations
(i--iii), $(u_1b)^2=(u_1a_4)^2=a'_3(a''_3)^{-1}$ and $a'_3h=ha'_3$
for $h\in\{a_1,a_4,b,u_1\}$. Since $\ker\rho_\ast$ is generated by
$a'_3(a''_3)^{-1}$, we obtain that
$$\rho_\ast(\M(F_{\alpha_3}))=\langle a_1,a_3,a_4,b,u_1\,|\,\mathrm{(i-ix)}\rangle.$$

Observe that $u_3$ preserves orientation of $\alpha_3$ and reverses
orientation of its neighborhood. It follows from sequence
(\ref{sequence}), that to obtain a presentation for
$\textrm{Stab}^+[\alpha_3]$ we have to add to the presentation for
$\rho_\ast(\M(F_{\alpha_3}))$ generator $u_3$ and relations
$u_3^2\in\rho_\ast(\M(F_{\alpha_3}))$ and
$u_3hu_3^{-1}\in\rho_\ast(\M(F_{\alpha_3}))$ for
$h\in\{a_1,a_3,a_4,b,u_1\}$. Thus
$$\textrm{Stab}^+[\alpha_3]=\langle a_1,a_3,a_4,b,u_1,u_3\,|\,\mathrm{(i-xv)}\rangle.$$
Analogously, since $t$ reverses the orientation of $\alpha_3$, we
obtain a presentation for $\textrm{Stab}[\alpha_3]$ by adding to the
above presentation generator $t$ and relations (xvi--xxii). \foorp

\begin{prop}\label{stab_b}
The stabilizer $S_{V_3}=\mathrm{Stab}[\beta]$ admits a presentation
with generators $a_1$, $a_2$, $a_3$, $b$, $t$, $w=u_1^{-1}u_3$, and
relations:

\noindent $\mathrm{(i)}\ ba_1=a_1b,\quad$ $\mathrm{(ii)}\
ba_2=a_2b,\quad$ $\mathrm{(iii)}\ ba_3=a_3b,\quad$ $\mathrm{(iv)}\
a_1a_3=a_3a_1,\quad$\\ $\mathrm{(v)}\ a_1a_2a_1=a_2a_1a_2,\quad$
$\mathrm{(vi)}\ a_2a_3a_2=a_3a_2a_3,\quad$ $\mathrm{(vii)}\
(a_1a_2a_3)^4=1,\quad$\\ $\mathrm{(viii)}\ t^2=1,\quad$
$\mathrm{(ix)}\ ta_1=a_1t,\quad$ $\mathrm{(x)}\ ta_2=a_2t,\quad$
$\mathrm{(xi)}\ ta_3=a_3t,\quad$\\ $\mathrm{(xii)}\
tbt=b^{-1},\quad$ $\mathrm{(xiii)}\ w^2=1,\quad$ $\mathrm{(xiv)}\
wa_1w=a_1^{-1},\quad$ $\mathrm{(xv)}\ wb=bw,$\\ $\mathrm{(xvi)}\
wa_3w=a_3^{-1},\quad$ $\mathrm{(xvii)}\
wa_2w=a_1a_3^{-1}a_2^{-1}a_3a_1^{-1},\quad$ $\mathrm{(xviii)}\
wt=tw.$

\noindent Relations \textrm{(i--xviii)} are consequences of
relations \textrm{(1--21)} in Theorem \ref{main}.
\end{prop}

\proof First we show that (i--xviii) are consequences of (1--21).
Relations (i--xii) appear among relations (1--21) in Theorem
\ref{main}; (xiii) follows from (13,14); (xiv) follows from (7) and
(i) in Proposition \ref{stab_d}; (xv) follows from (xiii) in
Proposition \ref{stab_a3}; (xvi) from (8) and (vi) in Proposition
\ref{stab_d}; (xviii) from (xiii), (13,18,19) and (xi) in
Proposition \ref{stab_d}. By relations (4,9)  we have:
$$wa_2w=u_1^{-1}u_3a_2u_3^{-1}u_1\stackrel{\small(9)}{=}u_1^{-1}a_2a_4^{-1}a_2^{-1}u_1
\stackrel{\small(4)}{=}u_1^{-1}a_4^{-1}a_2^{-1}a_4u_1.$$ From this
and (v) in Proposition \ref{stab_a3} we obtain that (xvii) is
equivalent to:
$$u_1a_2u_1^{-1}=a_4^{-1}a_1a_3^{-1}a_2a_3a_1^{-1}a_4.$$
From (5,15) we have
$$u_1a_2u_1^{-1}=(a_1a_2a_3)^{-2}u_3(a_1a_2a_3)^2a_2(a_1a_2a_3)^{-2}u_3^{-1}(a_1a_2a_3)^2,$$
and it is not difficult to check, that by (1,4)
$$(a_1a_2a_3)^2a_2(a_1a_2a_3)^{-2}=a_1a_3^{-1}a_2a_3a_1^{-1},$$
hence
$$u_1a_2u_1^{-1}=(a_1a_2a_3)^{-2}u_3a_1a_3^{-1}a_2a_3a_1^{-1}u_3^{-1}(a_1a_2a_3)^2\stackrel{\small(7,8,9)}{=}$$
$$=(a_1a_2a_3)^{-2}a_1a_3a_2a_4^{-1}a_2^{-1}a_3^{-1}a_1^{-1}(a_1a_2a_3)^2=$$
$$=a_3^{-1}a_2^{-1}a_1^{-1}\underline{a_3^{-1}a_2^{-1}a_3a_2}a_4^{-1}\underline{a_2^{-1}a_3^{-1}a_2a_3}a_1a_2a_3
\stackrel{\small(4)}{=}$$
$$=a_3^{-1}a_2^{-1}a_1^{-1}a_2\underline{a_3^{-1}a_4^{-1}a_3}a_2^{-1}a_1a_2a_3\stackrel{\small(1)}{=}
a_3^{-1}a_2^{-1}a_1^{-1}a_2a_4^{-1}a_2^{-1}a_1a_2a_3.$$ Thus (xvii)
is equivalent to:
$$a_3^{-1}a_2^{-1}a_1^{-1}a_2a_4^{-1}a_2^{-1}a_1a_2a_3=a_4^{-1}a_1a_3^{-1}a_2a_3a_1^{-1}a_4,$$
$$a_1^{-1}a_2a_4^{-1}a_2^{-1}a_1=a_2\underline{a_3a_4^{-1}a_1a_3^{-1}}a_2\underline{a_3a_1^{-1}a_4a_3^{-1}}a_2^{-1}
\stackrel{\small(1,2)}{\Leftrightarrow}$$
$$a_1^{-1}a_2a_4^{-1}a_2^{-1}a_1=a_2a_4^{-1}\underline{a_1a_2a_1^{-1}}a_4a_2^{-1}
\stackrel{\small(4)}{=}a_2a_4^{-1}a_2^{-1}a_1a_2a_4a_2^{-1}\stackrel{\small(9)}{\Leftrightarrow}$$
$$a_1^{-1}u_3a_2u_3^{-1}a_1=u_3a_2u_3^{-1}a_1u_3a_2^{-1}u_3^{-1}
\stackrel{\small(7)}{\Leftrightarrow}a_1^{-1}a_2a_1=a_2a_1a_2^{-1}\Leftarrow
(4).$$

The surface $F_{\beta}$ is torus with two holes. Let $b'$, $b''$
denote Dehn twists about its boundary components, such that
$\rho_\ast(b')=\rho_\ast(b'')=b$. Then, by the main theorem of
\cite{G}, $\M(F_{\beta})$ admits a presentation with generators
$a_1$, $a_2$, $a_3$, $b'$, $b''$ and relations (iv,v,vi),
$(a_1a_2a_3)^4=b'(b'')^{-1}$ and $b'h=hb'$ for
$h\in\{a_1,a_2,a_3\}$. Since $\ker\rho_\ast$ is generated by
$b'(b'')^{-1}$, we obtain that
$$\rho_\ast(\M(F_{\beta}))=\langle a_1,a_2,a_3,b\,|\,\mathrm{(i-vii)}\rangle.$$

Observe that $t$ preserves orientation of $\beta$ and reverses
orientation of its neighborhood. It follows from sequence
(\ref{sequence}), that to obtain a presentation for
$\textrm{Stab}^+[\beta]$ we have to add to the presentation for
$\rho_\ast(\M(F_{\beta}))$ generator $t$ and relations
(viii--xii). Then, since $w$ reverses the orientation of $\beta$,
we obtain a presentation for $\textrm{Stab}[\beta]$ by adding
generator $w$ and relations (xiii--xviii). \foorp

\begin{prop}\label{stab_e}
The stabilizer $S_{V_5}=\mathrm{Stab}[\varepsilon]$ is a subgroup of
$S_{V_3}$.
%
\end{prop}

\proof The surface $F_\varepsilon$ has two connected components. One
of them is torus with a hole, the other one is Klein bottle with a
hole containing $\beta$.
Let $h$ be any homeomorphism of $F$ which fixes $\varepsilon$. Then
$h$ fixes the connected components of $F_\varepsilon$. Since there
is only one isotopy class of unoriented non-separating two sided
curves in a Klein bottle with a hole, $h(\beta)$ and $\beta$ are
isotopic, hence $h\in\mathrm{Stab}[\beta]=S_{V_3}$.
\foorp

\begin{prop}\label{EG}
For $i\in\{1,\dots,11\}\backslash\{4,5\}$ the stabilizer $S_{E_i}$
is generated by the set $G_{E_i}$ defined in Table \ref{ch_edges}.
\end{prop}

\proof The surface $F_{(\alpha_3,\mu_1)}$ is projective plane with
three holes. By Theorem 7.5 of \cite{Sz1} and sequence
(\ref{sequence}), $\rho_\ast(\M(F_{(\alpha_3,\mu_1)}))$ is generated
by Dehn twists $a_3$, $a_4$, $y_1^{-1}a_4y_1$, $u_3^2$. Since $u_3$
preserves orientation of $\mu_1$ and $\alpha_3$ and reverses
orientation of a neighborhood of $\alpha_3$,
$\mathrm{Stab}^+[\alpha_3,\mu_1]$ is generated by $a_3$, $a_4$,
$y_1^{-1}a_4y_1$ and $u_3$. Since $t$ reverses orientation of both
$\alpha_3$ and $\mu_1$, while $y_1$ reverses orientation of $\mu_1$
only, $\mathrm{Stab}[\alpha_3,\mu_1]=S_{E_1}$ is generated by
$G_{E_1}=\{a_3, a_4, u_3, t, y_1\}$.

The surface $F_{(\alpha_3,\beta)}$ is a sphere with four holes. It
is a classical result (c.f. \cite{Bir}, Chapter 4) that the mapping
class group of such surface is generated by Dehn twists about the
boundary curves and three essential separating curves. In
$F_{(\alpha_3,\beta)}$ these essential curves may be taken as
$\alpha_1$, $(a_3^2a_2)^2(\alpha_1)$ and $\varepsilon$. Thus
$\rho_\ast(\M(F_{(\alpha_3,\beta)}))$ is generated by $a_3$, $b$ and
$a_1$, $(a_3^2a_2)^2a_1(a_3^2a_2)^{-2}$ and $e=(a_3^2a_2)^4$, by the
star relation (\ref{star1}). Suppose that
$h\in\Stab^+[\alpha_3,\beta]$ and $h$ reverses orientation of a
neighborhood of $\beta$. Then, since $F_\beta$ is orientable, $h$
also reverses orientation of a neighborhood of $\alpha_3$. Observe
that $(a_3^2a_2)^2t$ has this property. It follows that
$\Stab^+[\alpha_3,\beta]$ is generated by $b$, $a_3$, $a_1$, and
$(a_3^2a_2)^2t$, because
$(a_3^2a_2)^2a_1(a_3^2a_2)^{-2}=(a_3^2a_2)^2ta_1t^{-1}(a_3^2a_2)^{-2}$
and $(a_3^2a_2)^4=((a_3^2a_2)^2t)^2$, by relations (18,21) in
Theorem \ref{main}. Since $t$ preserves orientation of $\beta$ and
reverses orientation of $\alpha_3$, while $u_1^{-1}u_3$ reverses
orientation of $\beta$, $\Stab[\alpha_3,\beta]=S_{E_2}$ is generated
by $G_{E_2}=\{b,a_1,a_3, (a_3^2a_2)^2, t, u_1^{-1}u_3\}$.

The connected components of $F_{(\alpha_3,\delta)}$ are Klein bottle
with one hole and sphere with three holes. It is well known that the
mapping class group of a sphere with three holes is a free abelian
group of rank three generated by Dehn twists about the boundary
curves. It follows from sequence (\ref{sequence}) and Theorem A.7 of
\cite{S}, that $\rho_\ast(\M(F_{(\alpha_3,\delta)}))$ is generated
by $a_3$, $a_1$ and $u_1$. Observe that if
$h\in\Stab^+[\alpha_3,\delta]$ then $h$ fixes the components of
$F_\delta$, hence it preserves orientation of a neighborhood of
$\delta$. Since $u_3\in\Stab^+[\alpha_3,\delta]$ and it reverses
orientation of a neighborhood of $\alpha_3$,
$\Stab^+[\alpha_3,\delta]$ is generated by $a_3$, $a_1$, $u_1$ and
$u_3$. Suppose that $h\in\Stab[\alpha_3,\delta]$ and $h$ reverses
orientation of $\delta$. Then it induces an orientation reversing
homeomorphism of the orientable component of
$F_{(\alpha_3,\delta)}$, hence it reverses orientation of
$\alpha_3$. Since $t$ has this property,
$\Stab[\alpha_3,\delta]=S_{E_3}$ is generated by $G_{E_3}=\{a_3,
a_1, u_1, u_3, t\}$.

The surface $F_{(\mu_1,\varepsilon)}$ has two connected components.
One of the components is projective plane with two holes, hence its
mapping class group is free abelian group of rank two, generated by
Dehn twists abut its boundary components. The other component is
torus with one hole, hence its mapping class group is generated by
$a_2$ and $a_3$ (c.f \cite{G}). It follows from sequence
(\ref{sequence}) that $\rho_\ast(\M(F_{(\mu_1,\varepsilon)}))$ is
generated by $a_2$ and $a_3$. This group is equal to
$\mathrm{Stab}^+[\mu_1,\varepsilon]$ because every homeomorphism
fixing $\varepsilon$ must preserve its sides. Since $t$ reverses
orientation of $\mu_1$ and preserves orientation of $\varepsilon$,
while $u_3u_2u_3$ reverses orientation of $\varepsilon$,
$\mathrm{Stab}[\mu_1,\varepsilon]=S_{E_6}$ is generated by
$G_{E_6}=\{a_2, a_3, t, u_3u_2u_3\}$.

The surface $F_{(\mu_1,\delta)}$ has two connected components. One
of the components is projective plane with two holes, the other one
is Klein bottle with a hole. It follows from sequence
(\ref{sequence}) and Theorem A.7 of \cite{S}, that
$\rho_\ast(\M(F_{(\mu_1,\delta)}))$ is generated by $a_3$ and $u_3$.
Observe that any homeomorphism of $F$, which fixes $\mu_1$ and
$\delta$ must preserve the components of $F_\delta$. It follows that
if it preserves orientation of $\delta$, then it must also preserve
orientation of its neighborhood. Thus
$\rho_\ast(\M(F_{(\mu_1,\delta)}))=\mathrm{Stab}^+[\mu_1,\delta]$,
and $\mathrm{Stab}[\mu_1,\delta]=S_{E_7}$ is generated by
$G_{E_7}=\{a_3, u_3, t, y_1\}$.

The surface $F_{(\alpha_3,\alpha_1)}$ is sphere with four holes. 
Thus $\rho_\ast(\M(F_{(\alpha_3,\alpha_1)}))$ is generated by $a_1$,
$a_3$ and Dehn twists about curves $\delta$, $\beta$ and
$u_3(\beta)$, that is by $u_3^2$, $b$ and $u_3bu_3^{-1}$. Observe
that for $i\in\{1,3\}$, $u_i$ preserves orientation of $\alpha_i$
and reverses orientation of its neighborhood. Thus
$\mathrm{Stab}^+[\alpha_3,\alpha_1]$ is generated by $a_1$, $a_3$,
$b$, $u_1$ and $u_3$. Since $F_{(\alpha_3,\alpha_1)}$ is orientable,
any homeomorphism from $\mathrm{Stab}[\alpha_3,\alpha_1]$ which
reverses orientation of $\alpha_1$ must also reverse orientation of
$\alpha_3$. Observe that $t$ has this property, and thus
$\mathrm{Stab}[\alpha_3,\alpha_1]=S_{E_8}$ is generated by
$G_{E_8}=\{a_1, a_3, b, u_1, u_3, t\}$.

Both connected components of $F_{(\alpha_3,\alpha_4)}$ are
homeomorphic to the projective plane with two holes. It follows that
$\rho_\ast(\M(F_{(\alpha_3,\alpha_4)}))$ is generated by $a_3$ and
$a_4$. Note, that if $h\in\mathrm{Stab}^+[\alpha_3,\alpha_4]$
reverses orientation of a neighborhood of $\alpha_3$, then it must
interchange the components of $F_{(\alpha_3,\alpha_4)}$, and hance
also reverse orientation of a neighborhood of $\alpha_4$. Since
$u_3b$ has this property, it follows that
$\mathrm{Stab}^+[\alpha_3,\alpha_4]$ is generated by $a_3$, $a_4$
and $u_3b$. Observe that $u_1b$ reverses orientation of $\alpha_4$
and preserves orientation of $\alpha_3$, while $u_1t$ reverses
orientation of $\alpha_3$. Thus
$\mathrm{Stab}[\alpha_3,\alpha_4]=S_{E_9}$ is generated by
$G_{E_9}=\{a_3, a_4, u_3b, u_1b, u_1t\}$.

By Theorem 7.10 of \cite{Sz1}, $\M(F_{(\mu_1,\mu_2)})$ is generated
by $u_3$, $a_3$, $a_4$ and $y_2^2$. Observe that $y_2$ reverses
orientation of $\mu_2$ and preserves orientation $\mu_1$, while $t$
reverses orientation of $\mu_1$ and $\mu_2$. It follows that
$\mathrm{Stab}[\mu_1,\mu_2]=S_{E_{10}}$ is generated by
$G_{E_{10}}=\{u_3, a_3, a_4, t, y_2=u_2a_2\}$.

The surface $F_{(\mu_1,\mu_5)}$ is torus with two holes. Thus,
$\rho_\ast(\M(F_{(\mu_1,\mu_2)}))$ is generated by Dehn twists
$a_2$, $a_3$ and $a_4$ (c.f. \cite{G}). Since $F_{(\mu_1,\mu_5)}$ is
orientable, any homeomorphism from $\mathrm{Stab}[\mu_1,\mu_5]$
which reverses orientation of $\mu_1$ must also reverse orientation
of $\mu_5$. Observe that $u_3u_2u_3t$ has this property, and thus
$\mathrm{Stab}[\mu_1,\mu_5]=S_{E_{11}}$ is generated by
$G_{E_{11}}=\{a_2, a_3, a_4, u_3u_2u_3t\}$. \foorp

\section{\label{inj}Injectivity of $\Phi$.}
In this section we finish the proof of Theorem \ref{main} by showing
that the epimorphism $\Phi\colon\mathcal{G}\to\M(F)$ defined at the
end of Section \ref{preli} is injective.

For $i\in\{1,\dots,4\}$ let $\langle G_{V_i}|R_{V_i}\rangle$ be the
presentation for the stabilizer $S_{V_i}$ defined in Proposition
\ref{stab_mu1}, \ref{stab_d}, \ref{stab_a3}  or \ref{stab_b}, and
let $\langle G_{V_5}|R_{V_5}\rangle$ be any finite presentation for
$S_{V_5}$. For $j\in\{1,\dots,11\}\backslash\{4,5\}$ let $G_{E_j}$
be the generating set for $S_{E_j}$ defined in Table \ref{ch_edges},
and let $G_{E_4}$, $G_{E_5}$ be any finite generating sets for
$S_{E_4}$, $S_{E_5}$. For each $E\in X^1$ let
$G_{\overline{E}}=G_E$. Then $\M(F)$ admits the presentation defined
in Theorem \ref{presentation}. By Proposition \ref{stab_e},
$S_{V_5}\subset S_{V_3}$, hence each generator in $G_{V_5}$ may be
expressed in terms of $G_{V_3}$ and then the relations $R_{V_5}$
follow from $R_{V_3}$.
The relations
\begin{eqnarray}
&\label{RE1}&g_{E_i}=1=g_{\overline{E_i}}\quad\mathrm{for\ }i\le 7,\\
&\label{RE2}&g_{E_8}=(a_1a_2a_3)^2,\quad
g_{E_9}=a_2^{-1}u_2^{-1},\quad g_{E_{10}}=u_1,\quad
g_{E_{11}}=b^{-1}
\end{eqnarray}
obviously hold in $\M(F)$. It follows that the generating symbols
$g_E$ in relations $R^{(2)}$ and $R^{(3)}$ my be replaced by
expressions in generators $\bigcup_{i\le 4}G_{V_i}$. In order to
prove that $\Phi$ is injective it suffices to show that relations
$R_{V_i}$ for $i\le 4$, $R^{(2)}$ and $R^{(3)}$ are consequences of
relations (1--21) in Theorem \ref{main} and (\ref{RE1},\ref{RE2}).
For $R_{V_i}$ this is proved in Propositions \ref{stab_mu1},
\ref{stab_d}, \ref{stab_a3} and \ref{stab_b}. It remains to consider
$R^{(2)}$ and $R^{(3)}$.

\begin{figure}
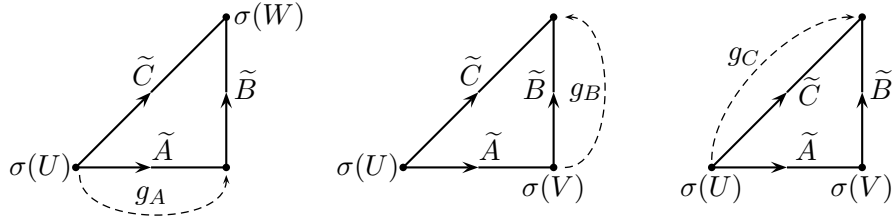

\begin{tabular}{ccc}
\input{tri3}&\input{tri4}&\input{tri5}
\end{tabular}
\caption{\label{tr1}Representatives of triangles with one loop.}
\end{figure}

\begin{prop}\label{tri}
For suitable choices of $\varphi$ and $\psi$, the relations
$R^{(3)}$ in Theorem \ref{presentation} corresponding to triangles
$T^\tau_i$ for $i<12$ are consequences of relations
(\ref{RE1},\ref{RE2}). The relation corresponding to $T_{12}$ is
equivalent to
\begin{equation}\label{RT12}
u_1u_2u_1=u_2u_1u_2
\end{equation}
and it is a consequence of relations $\mathrm{(1-21)}$ in Theorem
\ref{main}.
\end{prop}

\proof Let $T$ be a triangle in $X$ with edges $A$, $B$, $C$ and
vertices $U$, $V$, $W$.

{\it Case 1:} Suppose that $\widetilde{A}=\sigma(A)$,
$\widetilde{B}=\sigma(B)$, $\widetilde{C}=\sigma(C)$, $g_A=1$,
$g_B=1$, $g_C=1$. Then we can choose $\varphi=1$, $\psi=1$, so that
$\eta=1$ and the corresponding relation is $g_Ag_Bg_C^{-1}=1.$

{\it Case 2:} Suppose that $A$ is a loop, $\widetilde{A}=\sigma(A)$,
$\widetilde{C}=\sigma(C)=\sigma(B)$, $g_B=1$, $g_C=1$ and $g_A\in
S_W$. Then we can choose $\varphi=1$, $\psi=g_A^{-1}$, so that
$\eta=1$ and the corresponding relation is
$g_Ag_Bg_A^{-1}g_C^{-1}=1.$

{\it Case 3:} Suppose that $B$ is a loop,
$\widetilde{A}=\sigma(A)=\sigma(C)$, $\widetilde{B}=\sigma(B)$,
$g_A=1$, $g_C=1$ and $g_B\in S_U$. Then we can choose $\varphi=1$,
$\psi=1$, so that $\eta=g_B$ and the corresponding relation is
$g_Ag_Bg_C^{-1}=g_B.$

{\it Case 4:} Suppose that $C$ is a loop,
$\widetilde{A}=\sigma(A)=\sigma(\overline{B})$,
$\widetilde{C}=\sigma(C)$, $g_A=1$, $g_B=1$ and $g_C\in S_V$. Then
we can choose $\varphi=g_C$, $\psi=1$, so that $\eta=1$ and the
corresponding relation is $g_Ag_Cg_Bg_C^{-1}=1.$

Observe that for the representatives $\sigma(T_i^\tau)$ that we have
chosen in Section \ref{s_pres}, each of the $6$ triangles $T^\tau_i$
for $i\in\{3,5,7\}$ satisfies the assumptions of case 1. For
$i\notin\{3,5,7,12\}$, each of the $3$ triangles $T^\tau_i$
satisfies the assumptions of one of the cases 2, 3 or 4 (Figure
\ref{tr1}). It follows that the relations $R^{(3)}$ corresponding to
these triangles are consequences of (\ref{RE1},\ref{RE2}).

For triangle $T_{12}$ we have
$\widetilde{A}=[\mu_1,\mu_2]=\sigma(E_{10})=\sigma(A)=\sigma(B)=\sigma(C)$,
$\widetilde{B}=[\mu_2,\mu_3]$, $\widetilde{C}=[\mu_1,\mu_3]$,
$g_A=g_B=g_C=u_1$. We can take $\varphi=u_2$ and $\psi=u_2^{-1}$. We
claim that then $\eta=u_2$, so that the corresponding relation is
$g_Au_2g_Bu_2^{-1}g_C^{-1}=u_2$, which is equivalent to
(\ref{RT12}). Clearly it suffices to prove that (\ref{RT12}) is a
consequence of relations (1--21) in Theorem \ref{main}.

\noindent
$u_1u_2u_1\stackrel{\small(5,15,16)}{=}$\\
$$=(a_1a_2a_3)^{-2}u_3(a_1a_2a_3)^2a_3^{-1}a_2^{-1}u_3^{-1}a_2a_3(a_1a_2a_3)^{-2}u_3(a_1a_2a_3)^2
\stackrel{\small(7)}{=}
$$
$$=(a_1a_2a_3)^{-1}a_3^{-1}a_2^{-1}\underline{u_3a_2a_3u_3^{-1}}a_3^{-1}a_2^{-1}u_3a_2a_3a_1a_2a_3
\stackrel{\small(8,9)}{=}
$$
$$=(a_1a_2a_3)^{-1}a_3^{-1}a_4^{-1}a_2^{-1}a_3^{-2}a_2^{-1}u_3a_2a_3a_1a_2a_3.
$$
$u_2u_1u_2\stackrel{\small(5,15,16)}{=}$\\
$$=a_3^{-1}a_2^{-1}u_3^{-1}a_2a_3(a_1a_2a_3)^{-2}u_3(a_1a_2a_3)^2a_3^{-1}a_2^{-1}u_3^{-1}a_2a_3
\stackrel{\small(7)}{=}
$$
$$(a_1a_2a_3)^{-1}u_3^{-1}a_3^{-1}a_2^{-1}\underline{u_3a_2a_3u_3^{-1}}a_1a_2a_3
\stackrel{\small(8,9)}{=}$$
$$=(a_1a_2a_3)^{-1}u_3^{-1}a_3^{-1}a_4^{-1}a_2^{-1}a_3^{-1}a_1a_2a_3.$$
Now (\ref{RT12}) is equivalent to
$$a_3^{-1}a_4^{-1}a_2^{-1}a_3^{-2}a_2^{-1}u_3a_2a_3=\underline{u_3^{-1}a_3^{-1}a_4^{-1}}a_2^{-1}a_3^{-1}
\stackrel{\small(8,10)}{=}a_3a_4u_3a_2^{-1}a_3^{-1},$$
$$u_3a_2a_3^{2}a_2u_3^{-1}=a_2a_3^2a_2a_4a_3^2a_4
\stackrel{\small(8,9)}{\Leftrightarrow}
a_2a_4^{-1}a_2^{-1}a_3^{-2}a_2a_4^{-1}a_2^{-1}=a_2a_3^2a_2a_4a_3^2a_4,
$$
$$1=a_3^2a_2a_4a_3^2\underline{a_4a_2a_4a_2^{-1}}a_3^2a_2a_4
\stackrel{\small(4)}{=}(a_3^2a_2a_4)^3.
$$
It is not difficult to check that $(a_3^2a_2a_4)^3=1$ is a
consequence of (2,4,6).
 \foorp

\begin{prop}\label{R_edges}
The relations $R^{(2)}$ in Theorem \ref{presentation} corresponding
to edges of $X$ are consequences of (\ref{RE1},\ref{RE2}) and
relations $\mathrm{(1-21)}$ in Theorem \ref{main}.
\end{prop}

\proof For $i\in\{1,\dots,7\}$ we have
$g_{E_i}=g_{\overline{E_i}}=1$, thus relations corresponding to
$E_i$ identify, for each generator $g\in G_{E_i}$ of $S_{E_i}$, the
expression for $g$ in generators of $S_{i(E_i)}$ with the expression
in generators of $S_{t(E_i)}$. The relations corresponding to
$\overline{E_i}$ are the same, since $S_{E_i}=S_{i(E_i)}\cap
S_{t(E_i)}=S_{\overline{E_i}}$.
For $i\in\{8,\dots,11\}$ relations corresponding to the loop $E_i$
identify $g_{E_i}^{-1}gg_{E_i}$ as an element of $S_{i(E_i)}$ for
each $g\in G_{E_i}$.

Observe that all elements of $G_{E_1}$ except for $y_1$ appear as
generators in the presentations for $\Stab[\alpha_3]$ and
$\Stab[\mu_1]$. The only nontrivial relation corresponding to $E_1$
identifies expression for $y_1$ in generators of $\Stab[\alpha_3]$,
that is $u_1a_1$, with the expression in generators of
$\Stab[\mu_1]$ and it follows from (17):
$u_1a_1=u_2^{-1}u_3^{-1}ta_3^{-1}a_2^{-1}$.

The only nontrivial relation corresponding to $E_2$ identifies
$(a_3^2a_2)^2$ as an element of $\Stab[\alpha_3]$. By (17,21) in
Theorem \ref{main} we have $t=a_2a_3u_3u_2u_1a_1$, and
$$ta_1^{-1}u_1^{-1}
\stackrel{\small(16)}{=}
a_2a_3\underline{u_3a_3^{-1}a_2^{-1}u_3^{-1}a_2}a_3
\stackrel{\small(8,9)}{=}
a_2a_3^2a_2a_4a_3\stackrel{\small(2)}{=}a_2a_3\underline{a_3a_2a_3}a_4
\stackrel{\small(4)}{=}$$
$$=\underline{a_2a_3a_2}a_3a_2a_4
\stackrel{\small(4)}{=} a_3a_2a_3^2a_2a_4=a_3^{-1}(a_3^2a_2)^2a_4,$$
$$(a_3^2a_2)^2=a_3ta_1^{-1}u_1^{-1}a_4^{-1}\in\Stab[\alpha_3].$$

Note that all elements of $G_{E_3}$ appear as generating symbols for
$\Stab[\alpha_3]$ and $\Stab[\delta]$, so all relations
corresponding to $E_3$ are trivial.

Relations corresponding to $E_5$ identify the generators of
$\Stab[\varepsilon]$ as elements of $\Stab[\beta]$, because by
Proposition \ref{stab_e},
$\Stab[\beta,\varepsilon]=\Stab[\varepsilon]$.

Relations corresponding to $E_4$ are consequences of relations
corresponding to $E_{5}$ and $E_2$, because by Proposition
\ref{stab_e},
$\mathrm{Stab}[\alpha_3,\varepsilon]\subseteq\mathrm{Stab}[\alpha_3,\beta]$.

Relations corresponding to $E_6$ identify, for each $g\in G_{E_6}$,
the expression for $g$ in generators of $\mathrm{Stab}[\mu_1]$ with
the expression in generators of $\mathrm{Stab}[\varepsilon]$. But
every generator of $\mathrm{Stab}[\varepsilon]$ is identified with
an element of $\mathrm{Stab}[\beta]$, by relations corresponding to
$E_{5}$. The only nontrivial relation identifies $u_3u_2u_3$ as an
element of $\mathrm{Stab}[\beta]$. By (17,21) we have
$t=a_1a_2a_3u_3u_2u_1$, and
$$u_3u_2u_3=(a_1a_2a_3)^{-1}tu_1^{-1}u_3\in\mathrm{Stab}[\beta].$$

The only nontrivial relation corresponding to $E_7$ identifies
expression for $y_1$ in generators of $\mathrm{Stab}[\delta]$, that
is $u_1a_1$, with an expression in generators of
$\mathrm{Stab}[\mu_1]$. Such relation can be derived from (17).

Relations corresponding to $E_8$ are: $s^{-1}a_1s=a_3$,
$s^{-1}a_3s=a_1$ $s^{-1}bs=b$, $s^{-1}u_1s=u_3$, $s^{-1}u_3s=u_1$,
$s^{-1}ts=t$, where $s=g_{E_8}=(a_1a_2a_3)^2$, and they all follow
from relations in Proposition \ref{stab_d} and (3) in Theorem
\ref{main}.

Relations corresponding to $E_9$ are
$u_2a_2ga_2^{-1}u_2^{-1}\in\mathrm{Stab}[\alpha_3],$ for $g\in
G_{E_9}$. It can be checked that
$u_2a_2a_4a_2^{-1}u_2^{-1}=a_3^{-1}$ and
$u_2a_2a_3a_2^{-1}u_2^{-1}=ta_4^{-1}t$. Observe that the last two
relations involve only generators from $\mathrm{Stab}[\mu_1]$, and
hence they are consequences of relations in Proposition
\ref{stab_mu1}.

From (3,11,16) we have
\begin{equation}\label{u_2b}
(u_2^{-1}b)^2=1.
\end{equation}

Using (xiii) in Proposition \ref{stab_mu1}, (\ref{u_2b}) and (3), we
obtain:
$$u_2a_2u_3ba_2^{-1}u_2^{-1}
=a_2^{-1}u_2u_3\underline{bu_2^{-1}}a_2
\stackrel{\small(\ref{u_2b})}{=}a_2^{-1}u_2u_3u_2b^{-1}a_2=a_2^{-1}u_2u_3u_2a_2b^{-1}.$$
By relations in Theorem \ref{main} we have:
$$a_2^{-1}u_2u_3u_2a_2\stackrel{\small(16)}{=}
\underline{a_2^{-1}a_3^{-1}a_2^{-1}}u_3^{-1}a_2a_3\underline{u_3a_3^{-1}}a_2^{-1}u_3^{-1}a_2a_3a_2
\stackrel{\small(4,8)}{=}$$
$$a_3^{-1}a_2^{-1}\underline{a_3^{-1}u_3^{-1}}a_2a_3^2\underline{u_3a_2^{-1}u_3^{-1}a_2}a_3a_2
\stackrel{\small(8,9)}{=}
a_3^{-1}u_3^{-1}\underline{u_3a_2^{-1}u_3^{-1}}a_3a_2a_3^2a_2a_4a_3a_2$$
$$\stackrel{\small(9)}{=}
a_3^{-1}u_3^{-1}a_2a_4\underline{a_2^{-1}a_3a_2a_3^2}a_2a_4a_3a_2
\stackrel{\small(4)}{=}
a_3^{-1}u_3^{-1}a_2a_4a_3\underline{a_2a_3a_2a_4a_3}a_2$$
$$\stackrel{\small(4)}{=}
a_3^{-1}u_3^{-1}a_2a_4a_3^2a_2\underline{a_3a_4a_3}a_2
\stackrel{\small(2)}{=}
a_3^{-1}u_3^{-1}(a_2a_4a_3^2)^3a_3^{-2}a_4^{-1}.$$ It is not
difficult to check, that by (2,4,6),
$(a_2a_4a_3^2)^3=(a_4a_2a_3)^4=1$. Thus
$$u_2a_2u_3ba_2^{-1}u_2^{-1}=\underline{a_3^{-1}u_3^{-1}}a_3^{-2}a_4^{-1}b^{-1}
\stackrel{\small(8)}{=}(ba_4a_3u_3)^{-1}\in\mathrm{Stab}[\alpha_3].$$

Before we describe the remaining two relations (for $g=u_1b,u_1t$)
we will show that relation
\begin{equation}\label{local1}
a_2^{-1}u_1u_2u_1a_2=a_3wt,
\end{equation} where $w=u_1u_3^{-1}$, is a consequence of
relations in Theorem \ref{main}. By (17,21) we have
$t=a_1a_2a_3u_3u_2u_1$, and (\ref{local1}) is equivalent to
$$a_3w=a_2^{-1}u_1u_2u_1t^{-1}a_2=a_2^{-1}u_1u_3^{-1}a_3^{-1}a_2^{-1}a_1^{-1}a_2=
a_2^{-1}wa_3^{-1}a_2^{-1}a_1^{-1}a_2,$$
$$wa_3w=wa_2^{-1}wa_3^{-1}a_2^{-1}a_1^{-1}a_2.$$
By (xvi,xvii) in Proposition \ref{stab_b}, this is equivalent to
$$a_3^{-1}=a_1a_3^{-1}a_2a_3a_1^{-1}a_3^{-1}a_2^{-1}a_1^{-1}a_2,$$
and it is easy to check, that the last relation is a consequence of
(1,4).

Now, from (xiii) in Proposition \ref{stab_mu1}, (\ref{u_2b}) and
(\ref{RT12}), we obtain:
$$u_2a_2u_1ba_2^{-1}u_2^{-1}=a_2^{-1}u_2u_1bu_2^{-1}a_2=a_2^{-1}u_2u_1u_2b^{-1}a_2=a_2^{-1}u_1u_2u_1a_2b^{-1},$$
hence, by (\ref{local1})
$$u_2a_2u_1ba_2^{-1}u_2^{-1}=a_3wtb^{-1}\in\mathrm{Stab}[\alpha_3].$$
Similarly, using (xiv) in Proposition \ref{stab_mu1}, we have
$$u_2a_2u_1ta_2^{-1}u_2^{-1}=a_2^{-1}u_2u_1u_2a_2t=a_2^{-1}u_1u_2u_1a_2t=a_3w\in\mathrm{Stab}[\alpha_3].$$

The relations corresponding to $E_{10}$ are $u_1^{-1}u_3u_1=u_3$,
$u_1^{-1}a_3u_1=a_3$, $u_1^{-1}a_4u_1=u_3^{-1}a_4u_3$,
$u_1^{-1}tu_1=tu_3^2$ and
$u_1^{-1}u_2a_2u_1\in\mathrm{Stab}[\mu_1]$. First four relations are
easy consequences of relations in Proposition \ref{stab_a3}. By
(4,8,16) in Theorem \ref{main}
$$u_2a_2\stackrel{\small(16)}{=}a_3^{-1}a_2^{-1}u_3^{-1}\underline{a_2a_3a_2}
\stackrel{\small(4)}{=}a_3^{-1}a_2^{-1}\underline{u_3^{-1}a_3}a_2a_3
\stackrel{\small(8)}{=}a_3^{-1}a_2^{-1}a_3^{-1}u_3^{-1}a_2a_3,$$
thus by (xvi,xvii) in Proposition \ref{stab_b} and (13)
$$wu_2a_2w=a_3a_1a_3^{-1}a_2a_3a_1^{-1}a_3u_3^{-1}a_1a_3^{-1}a_2^{-1}a_3a_1^{-1}a_3^{-1},$$
which is equivalent, by (1,7), to
$$wu_2a_2w=(a_1a_2a_3)a_3u_3^{-1}(a_1a_2a_3)^{-1}.$$
By (xiii,xiv,xv) in Proposition \ref{stab_d} we have
$$(a_1a_2a_3)a_3u_3^{-1}(a_1a_2a_3)^{-1}=(a_1a_2a_3)^{-1}a_1u_1^{-1}(a_1a_2a_3),$$ hence, using (i) in
Proposition \ref{stab_d},
$wu_2a_2w=a_3^{-1}a_2^{-1}(u_1a_1)^{-1}a_2a_3$ and
$$u_1^{-1}u_2a_2u_1=u_3^{-1}a_3^{-1}a_2^{-1}(u_1a_1)^{-1}a_2a_3u_3^{-1}.$$
It remains to notce that $u_1a_1$ may be written in generators of
$\mathrm{Stab}[\mu_1]$ using (17):
$u_1a_1=u_2^{-1}u_3^{-1}ta_3^{-1}a_2^{-1}$.

The relations corresponding to $E_{11}$ are $ba_2b^{-1}=a_2$,
$ba_3b^{-1}=a_3$, $ba_4b^{-1}=u_3^{-1}a_4u_3$ and
$bu_3u_2u_3tb^{-1}=(u_3u_2u_3)^{-1}t$. First two follow from (3),
third follows from (iii,xiv) in Proposition \ref{stab_a3}, fourth
follows from (11,20) and (\ref{u_2b}):
$$bu_3u_2u_3tb^{-1}\stackrel{\small{(11,20)}}{=}
u_3^{-1}b^{-1}u_2b^{-1}u_3^{-1}t\stackrel{\small(\ref{u_2b})}{=}(u_3u_2u_3)^{-1}t.$$\foorp

\end{document}